\documentclass{article}
\usepackage{amsfonts}
\usepackage{amscd}
\usepackage{amsmath}
\usepackage{verbatim}
\usepackage{amssymb}
\usepackage{epic}

\newtheorem{theorem}{Theorem}
\newtheorem{theorema}{Theorem}
\newtheorem{theoremb}{Theorem}
\newtheorem{theoremc}{Theorem}

\newtheorem{cor}[theorem]{Corollary}

\newtheorem{dfn}[theorema]{Definition}
\newtheorem{examp}[theoremc]{Example}

\newtheorem{prop}[theorem]{Proposition}
\newtheorem{remark}[theoremb]{Remark}

\newenvironment{proof}[1][Proof]{\textbf{#1.} }{\qed}

\renewcommand\a{\alpha}

\renewcommand\c{\circ}
\newcommand\C{{\cal C}}
\newcommand\CC{{\mathbb C}}

\renewcommand\d{\delta}
\newcommand\D{{\cal D}}
\newcommand\e{\eta}
\newcommand\E{{\cal E}}

\newcommand\Fg{{\frak F}}

\newcommand\g{{\frak g}}

\newcommand\hps{\hskip-16pt . \hskip2pt}
\newcommand\hpss{\hskip-13.5pt . \hskip2pt}

\renewcommand\ll{\lambda}
\newcommand\La{\Lambda}

\newcommand\oo{\omega}

\newcommand\op[1]{\mathop{\rm #1}\nolimits}
\newcommand\ot{\otimes}
\newcommand\p{\partial}

\newcommand\po{$\!\!\!{\bf .}$ }
\renewcommand\P{\Phi}

\newcommand\R{{\mathbb R}}
\newcommand\St{\op{St}^l_a}
\renewcommand\t{\tau}
\newcommand\tg{{\frak t}}
\newcommand\tim-
\newcommand\te{\theta}

\newcommand\T{\Theta}
\newcommand\ups{\upsilon}

\newcommand\uu{{\cal U}_G}

\newcommand\vg{{\frak v}}
\newcommand\vp{\varphi}

\newcommand\we{\wedge}
\newcommand\x{\xi}

\newcommand\z{\sigma}

\def\Rom#1{\uppercase\expandafter{\romannumeral#1}}

\newcommand\1{{\bf 1}}

\newcommand\qed{\phantom{\underline{y}}\hfill\hfill$\square$}

\newcommand\bib[1]{\bibitem[#1]{#1}}
\newcommand\abz{\hspace{12pt}}

\begin{document}

\title{\bf Invariants of pseudogroup actions:\\
Homological methods and \\
Finiteness theorem}
\author{\textsc{Boris Kruglikov and Valentin Lychagin}
\medskip \\
\textsf{Mat-Stat. Dept., University of Troms\o, Norway}}

\maketitle

\begin{abstract}
We study the equivalence problem of submanifolds with respect to a
transitive pseudogroup action. The corresponding differential
invariants are determined via formal theory and lead to the
notions of $l$-variants and $l$-covariants, even in the case of
non-integrable pseudogroup. Their calculation is based on the
cohomological machinery: We introduce a complex for covariants,
define their cohomology and prove the finiteness theorem. This
implies the well-known Lie-Tresse theorem about differential
invariants. We also generalize this theorem to the case of
pseudogroup action on differential equations.
 \footnote{
pseudogroup, differential invariants, Tresse derivative,
equivalence, Lie equation, Spencer cohomology.}%
\end{abstract}

\section*{Introduction}\label{sec0}

 \abz
Transformation groups were introduced by S.\,Lie \cite{Li1} in his
study of invariants of differential equations. More formal and
general notion of pseudogroup was developed by E.\,Cartan
\cite{C1}. These infinite groups of Lie and Cartan (see also
\cite{H,Tr1,V}) play nowadays a central role in geometry and
analysis.

A pseudogroup $G\subset\op{Diff}_{\op{loc}}(M)$ acting on a
manifold $M$ consists of a collection of local diffeomorphisms
$\vp$, each bearing own domain of definition $\op{dom}(\vp)$ and
range $\op{im}(\vp)$, that satisfies the following properties:
 \begin{enumerate}
  \item $\op{id}_M\in G$ and $\op{dom}(\op{id}_M)=\op{im}(\op{id}_M)=M$,
  \item If $\vp,\psi\in G$, then $\vp\circ\psi\in G$ whenever
$\op{dom}(\vp)\subset\op{im}(\psi)$,
  \item If $\vp\in G$, then $\vp^{-1}\in G$ and $\op{dom}(\vp^{-1})=\op{im}(\vp)$,
  \item $\vp\in G$ iff for every open subset $U\in\op{dom}(\vp)$
the restriction $\left.\vp\right|_U\in G$,
  \item The pseudogroup is of order $l$ if this is the minimal number such
that $\vp\in G$ whenever for each point $a\in\op{dom}(\vp)$ the
$l$-jet is admissible: $[\vp]^l_a\in G^l$.
 \end{enumerate}

The latter property means that a pseudogroup is defined by
differential equations of maximal order $l$ and will be explained
below. It uses the infinitesimal language. In fact, for most
infinite pseudogroups the only comprehensive way to introduce the
notion of continuity is via the prolongation technique.

A transformation $\vp\in G$ defines a map ($l$-th prolongation) of
the space of jets of codimension $r$ submanifolds
$\vp_{(l)}:J^l_r(M)\to J^l_r(M)$, which obeys the following
property:
 $$
(\vp\circ\psi^{-1})_{(l)}=\vp_{(l)}\circ\psi_{(l)}^{-1}.
 $$
This property is fundamental and being coupled with the formal
theory of differential equations leads to a deep understanding of
pseudogroups, cf.\ \cite{E,Lib,SS}.

In this paper we develop a more general notion of infinitesimal
pseudogroup. First of all, we don't require a pseudogroup to be
integrable, and so consider an $l$-pseudogroup as a transformation
group in finite jets. Even such finite order pseudogroups are
important in producing invariants for differential equations and
curvatures for geometric structures.

Next, we consider pseudogroups acting on differential equations
and this relates the theory with the original approach of S. Lie,
which is now called the theory of symmetries and is extensively
used in explicit integration of differential equation. Using this
technique the known invariant differential operators, connections
etc can be obtained.

Finally we do no require that pseudogroups act effectively. In our
approach the stationary sub-pseudogroups appear naturally. This is
convenient for applications, since we can consider then general
representations.

Formal integrability criterion for infinitesimal pseudogroups is
based on the well-developed algebraic machinery, described in the
paper. The passage from formal integrability to the local one is
not automatic and is generically wrong. However the former implies
the latter in the following cases:
 \begin{itemize}
  \item[-]
Finite type pseudogroups (the symbol $\mathfrak{g}^k\equiv0$ for
big $k$). This condition implies that the integrated pseudogroup
is finite-dimensional (Lie group).
 \item[-]
Analytic pseudogroups. It is a consequence of Cartan-K\"ahler
theorem, which holds for general differential equations
\cite{M,KLV}.
 \item[-]
Elliptic pseudogroups of analytic type, see \cite{S,M}.
 \item[-]
Transitive flat pseudogroups, see \cite{BM,P}.
 \end{itemize}
Only in special cases the global integrability (or equivalence)
problem can be handled, see \cite{S,GS2,Ta}.

Like the classical invariant theory, the theory of invariants of
pseudogroup actions exists. For Lie pseudogroups this is the
well-known differential invariants theory. Generally the number of
these differential invariants is infinite (the opposite case is of
much importance, was discussed in our preceding paper \cite{KL2}
and will be reviewed and continued here). But the whole algebra
can be finitely generated (under certain regularity assumptions).
This can be done via Lie approach with a finite number of
invariant differentiations generating all invariants from a finite
number of them (\cite{Li3,Li4})) or with Tresse method of
differentiating some invariants by the others (\cite{Tr1}).

In this paper we address the question of calculation of the
differential invariants and estimation of their number. We develop
the theory of pseudogroups actions on submanifolds, which are
subject to certain differential relations. In other words, we
consider the action of a pseudogroup on a system of differential
equations, which we naturally identify with a submanifold in the
space of jets.

On the level of finite jets we are lead to what we call
$l$-covariants. Their calculus is governed by certain cohomology
theory similar to the formal Spencer cohomology. We exploit this
in relation to the problem of equivalence of submanifolds under
the action.

Our approach gives the finiteness theorem for the cohomology of
covariants and this in turn implies Lie (\cite{Li2,Li5}) and
Tresse (\cite{Tr1}) theorems (proved by Kumpera \cite{Kum}, see
also modifications in \cite{Ov,Ol,MMR}) and their generalization
for the action of pseudogroups on differential equations (note
that in Kumpera's approach the pseudogroup is integrable and he
considers the integrated sheaf of vector fields, while we start
with general pseudogroup and calculate precisely how many
invariants we need on each jet level).

Lie pseudogroups consist of pseudo-automorphisms of geometric
structures. We provide a series of calculations for them. An
important class form the symmetries of differential equations,
realized as transformations preserving the Cartan distribution on
the equation. They are very useful for integration of the given
equation.

\section{\hps Pseudogroups and equivalence}\label{S1}

 \abz
In this section we discuss the general introduction of
pseudogroups, developing the ideas of \cite{GS2,Kur,S,Ta}. This
approach was proposed in \cite{KL2}.

\subsection{\hpss Formal pseudogroups}\label{S11}

 \abz
Let $M$ be a smooth manifold and $J^l_r(M)$ be the corresponding
jet space. Its points $a_l$ are the $l$-jets $[N]^l_a$ of
submanifolds $N\subset M$ of codimension $r$ at $a\in N$.

Denote the natural projections by $\rho_{i,j}:J^i_r(M)\to
J^j_r(M)$, $i\ge j$. It $j\ge1$ the fibers bear a canonical affine
structure (\cite{KLV,Ly}), associated with the vector structure,
described below. It is sufficient to specify it for
$\Fg(a_{l-1})=\rho_{l,l-1}^{-1}(a_{l-1})$.

Denote $\tg_a=T_aN=[N]^1_a$ and $\vg_a=T_aM/T_aN$. Let $a_l\in
J^l_r(M)$, $a_{l-1}=\rho_{l,l-1}(a_l)$. Then
$T_{a_l}\Fg(a_{l-1})\simeq S^l\tg^*_a\ot\vg_a$ and we get the
exact sequence:
 $$
0\to S^l\tg^*_a\ot\vg_a\to
T_{a_l}J^l_r(M)\stackrel{(\rho_{l,l-1})_*} \longrightarrow
T_{a_{l-1}}J^{l-1}_r(M)\to 0.
 $$

For a vector bundle $\rho:E\to B$ of rank $r$, the corresponding
space of jets of sections $J^l\rho$ is an open subset in
$J^l_r(E)$. In particular, we realize the jet space for maps
$J^l(N,M)$. Denote by $D^l(M)\subset J^l(M,M)$ the open dense
subset, consisting of the $l$-jets of local diffeomorphisms. Being
equipped with the partially defined composition operation, it is
an example of finite order pseudogroup.

To define this notion in general, recall some basic facts from the
geometric theory of differential equations, see \cite{KLV,Ly,Gu}
and Appendix \ref{A1} for details. The prolongation of
differential equation $\E\subset J^l_r(M)$ is defined as
 \begin{multline*}
\E^{(1)}=\{a_{l+1}=[N]_a^{l+1}\,|\,\text{ for }N\subset M\text{ if
codimension }r\text{ the jet-extension }\\
j_l(N)\subset J^l_r(M)\text{ is tangent to }\E\text{ at
}a_l\}\subset J^{l+1}_r(M).
 \end{multline*}
This can be equivalently written as
$\E^{(1)}=\{a_{l+1}\,|\,L(a_{l+1})\subset T_{a_l}\E\}$, where for
$a_{l+1}=[N]_a^{l+1}$ we set: $L(a_{l+1})=T_{a_l}j_l(N)$,
$a_l=\rho_{l+1,l}(a_{l+1})$.

The other prolongations are defined inductively:
$\E^{(i)}=(\E^{(i-1)})^{(1)}$.

 \begin{dfn}\po\label{pseudogr}
An $l$-pseudogroup is a collection of (non-empty) subbundles
$G^j\subset D^j(M)$, $0<j\le l$, such that the following
properties are satisfied:
 \begin{enumerate}
  \item If $\vp_j,\psi_j\in G^j$, then
$\vp_j\circ\psi_j^{-1}\in G^j$ whenever defined,
  \item The map $\rho_{j,j-1}:G^j\to G^{j-1}$ is a bundle for every $j\le l$.
 \end{enumerate}
As usual in the differential equations theory we assume
$G^0=D^0(M)= M\times M$, which is equivalent to {\em
transitivity\/} of the pseudogroup action.

An $l$-pseudogroup is called $l$-integrable if
$G^j\subset(G^{j-1})^{(1)}$ for all $0<j\le l$.
 \end{dfn}

Note that assumption 1 implies that $\op{id}^j_M\in G^j$ and
$\vp_j\in G^j\Rightarrow\vp_j^{-1}\in G^j$.

Pseudogroups $G=\{G^j\}_{j=1}^l$ defined by this approach can be
studied for integrability by the standard prolongation-projection
method, see \cite{GS1,GS2,Kur,KLV,S,Ta} and Appendix \ref{A1}.

Denote $G^j_{a,b}=\{\vp_j\in G^j\,|\,\vp_0(a)=b\}$,
$G^j_a=G^j_{a,a}$ -- the subgroup of $G^j$ and
$\mathfrak{G}^j_a=\op{Ker}[\rho_{j,j-1}:G^j_a\to G^{j-1}_a]$ --
its (normal) subgroup, which is abelian for $j>1$ and for $j=1$:
$\mathfrak{G}^1_a=G^1_a\subset\op{Gl}(T_aM)$.

 \begin{dfn}\po
Let $\vp_j\in G^j$ be a point and $\rho_{j,0}(\vp_j)=(a,b)\in
M\times M$. The symbol of the pseudogroup $G$ is given by:
 $$
\g^j(\vp_j)=\op{Ker}\bigl[(\rho_{j,j-1})_*:T_{\vp_j}G^j\to
T_{\vp_{j-1}}G^{j-1}\bigr].
 $$
It can be viewed as a subspace $\g^j(\vp_j)\subset S^j(T^*_aM)\ot
T_bM\stackrel{1\ot\vp_1^{\!-1}\!\!\!}\simeq S^j(T^*_aM)\ot T_aM$,
and in the last form is identified with the Lie algebra $\g^j_a$
of the Lie group $\mathfrak{G}^j_a$.
 \end{dfn}

An $l$-pseudogroup $G$ is called {\em formally integrable\/} if it
is $l$-integrable, for every $j>l$ the prolongation
$G^j=(G^l)^{(j-l)}$ exists, are $j$-pseudogroups and the
projections $\rho_{j,j-1}:G^j\to G^{j-1}$ are vector bundles.

Similar to the differential equations theory (\cite{Go,Gu,S}), a
criterion of formal integrability can be formulated in terms of
the Spencer $\d$-complex:
  \begin{equation}\label{hhh1}
0 \to\g^l_a\stackrel\d\to\g^{l-1}_a\ot T^*_aM\stackrel\d\to
\dots\stackrel\d\to\g^{l-j}_a\ot\La^jT^*_aM \stackrel\d\to \dots
  \end{equation}
Its bi-graded cohomology groups are denoted by $H^{l-j,j}(G)$ or
$H^{l-j,j}(\g)$. Poincar\'e $\d$-lemma states that $\dim
H^{*,*}(G)<\infty$, so that all cohomology groups $H^{i,j}(\g)$
eventually vanish. But some groups are non-zero.

The obstructions to formal integrability of the $l$-pseudogroup
$G$, considered as a differential equation, are some elements
$W_j(G)\in H^{j-1,2}(G)$, called Weyl tensors (or curvatures),
defined via the jet-spaces geometry. We do not need their precise
form here and so refer for the definition to \cite{Ly}.

 \begin{theorem}\po
Let $G$ be an $l$-pseudogroup. Suppose the symbols $\g^j$ over
$G^l$ form a vector bundle and all the Weyl tensors $W_j$ vanish
identically for all $j\ge l$. Then the pseudogroup is formally
integrable.
 \end{theorem}

 \begin{proof}
The hypotheses imply integrability of $G$ as a differential
equation, see \cite{Ly}. We need to check that the obtained system
$\{G^j\}_{j=0}^\infty$ is a pseudogroup, i.e.\ to check all the
requirements of definition 1.

Let $G^{j+1}=(G^j)^{(1)}$. Obviously, the unit is in $G^{j+1}$.
Let $\vp_{j+1}\in G^{j+1}_{a,b}$, $\psi_{j+1}\in G^{j+1}_{b,c}$
and $\chi_{j+1}=\psi_{j+1}\circ \vp_{j+1}$. We need to show that
$\chi_{j+1}\in G^{j+1}$. This is equivalent to
$L(\chi_{j+1})\subset T_{\chi_j}G^j$.

To prove the inclusion consider the multiplication operator
$m_j:G^j\times G^j\to G^j$. It has the differential:
 $$
T_{\psi_j}G^j\oplus T_{\vp_j}G^j\stackrel{dm_j}\longrightarrow
T_{\chi_j}G^j.
 $$
The two summands on the left contain the subspaces $L(\psi_{j+1})$
and $L(\vp_{j+1})$ respectively. But
 $$
L(\psi_{j+1})\oplus L(\vp_{j+1})\stackrel{dm_j}\longrightarrow
L(\psi_{j+1}\vp_{j+1})
 $$
for any $\vp_{j+1},\psi_{j+1}\in D^{j+1}(M)$ such that the
composition is defined. Moreover the multiplication operator with
one fixed argument is invertible. So the above arguments yield
that $\vp_{j+1}\in G^{j+1}$ implies $\vp_{j+1}^{-1}\in G^{j+1}$.
The claim follows.
 \end{proof}
\vspace{4pt}

An $l$-pseudogroup $G$ is called $q${\em -acyclic\/} if
$H^{i,j}(G)=0$ for $i\ge l$, $0\le j\le q$. An $\infty$-acyclic
pseudogroup is called {\em involutive\/}.
For such pseudogroups $G$ investigation of formal integrability
involves only one obstruction $W_l(G)$.

If a pseudogroup $G$ is formally integrable we obtain its infinite
prolongation
 $$
G^\infty=\op{lim}_{\op{proj}}(G^l,\rho_{l,l-1}),
 $$
which is called infinitesimal or {\em formal pseudogroup\/}. If
there is local integrability (smooth or analytic), as described in
the introduction, we refer to the pseudogroup as to {\em
integrable\/}.

Lie pseudogroups are characterized by the property that they can
be restored from the corresponding equation on jets of vector
fields, precisely in the way the Lie groups and algebras are
related, see \cite{KS}. The corresponding Lie equations for such
pseudogroups are always linear.

 \begin{examp}\po\rm
The group of complex fractional-linear transformations of $S^2=\CC
P^1$ (or real transformations of $S^1=\R P^1$) is an integrable
pseudogroup of finite type and order 3. In fact its Lie algebra is
represented as the algebra of quadratic-polynomial vector fields
on the line: $\g=\{\x=(c_0+c_1z+c_2z^2)\p_z\}$.
 \end{examp}

 \begin{examp}\po\rm
Consider the pseudogroup of local plane transformations
 $$
(q,p)\mapsto(F(q),p/F'(q)).
 $$
This is a Lie pseudogroup of infinite type and order 1. Indeed, it
consists of transformations from $T^*\R^1\simeq\R^2(q,p)$
preserving the Liouville form $p\,dq$. The generating field has
the form $\xi=f(q)\p_q-f'(q)p\,\p_p$. If we vary $F(q)$ in a
finite-dimensional subgroup of $\op{Diff}(\R^1)$, the pseudogroup
becomes of finite type.
 \end{examp}

 \begin{examp}\po\rm
Let $\E$ be a geometric structure (\cite{Gu,Ly}) and $G$ be its
Lie pseudogroup of the jets-automorphisms. If the structure $\E$
is integrable (flat), the pseudogroup is integrable as well. It
can be of finite or infinite type depending on the geometric
structure (\cite{Ko}). It has the same order as the structure
$\E$. When the geometric structure is non-integrable, the order of
the pseudogroup $G$ can increase and it can readily be
non-integrable (formally or locally).
 \end{examp}

One of the most important Lie pseudogroups consist of Lie
transformations on the jet-space $M=J^k\pi$ of some bundle
$\pi:E_\pi\to B$ (\cite{KLV}). It has order 1 and infinite type.
We will discuss this example in detail in \S\ref{S4} and Appendix
\ref{A2}.

\subsection{\hpss Pseudogroup action}\label{S12}

 \abz
A pseudogroup $G$ is represented by the action on local
submanifolds $N\subset M$ of codimension $r$. A formal pseudogroup
acts on the space $J^l_r(M)$.

The equivalence problem is to realize when a submanifold
$N_1\subset M$ can be transformed to a submanifold $N_2\subset M$
by a map $\vp\in G$. For formal pseudogroups we consider the
infinitesimal problem for $l$-jets and $l$-pseudogroups:

 \begin{dfn}\po
We say that $l$-jets of two submanifolds $N_1$ and $N_2$ at the
points $a,b\in M$ are $G$-equivalent if $\vp_l[N_1]^l_a=[N_2]^l_b$
for some $\vp_l\in G^l_{a,b}$.
 \end{dfn}

For transitive pseudogroups the equivalence problem reduces to the
case $a=b$. We assume this and begin subsequently equalizing the
jets of submanifolds.

The pseudogroup $D^l(M)$ and hence $G^l$ act on the space
$J^l_r(M)$ by the formula
$\vp_{(l)}:[N]_a^l\mapsto[\vp(N)]_{\vp(a)}^l$. These actions obey
the relation: $\rho_{l,s}\circ\vp_{(l)}=\vp_{(s)}\circ\rho_{l,s}$.

Consequently, the group $\mathfrak{G}^l_a$ acts on $\Fg(a_{l-1})$.
For $l=1$ this action is generated by the linear collineations in
the Grassmannians. The action is affine for $l>1$:
 $$
f\mapsto\ll(\te)+f,\quad f\in T_{a_l}\Fg(a_{l-1}).
 $$
Here $\ll$ is the induced linear representation of the Lie algebra
$\g^l_a$, which is naturally the restriction-factorization map:
 \begin{equation}\label{lambda}
\ll:\ S^lT_a^*M\ot T_aM\to T_{a_k}\Fg(a_{l-1}),\quad
\te\mapsto\bar\te\in S^l\tg^*_a\ot\vg_a.
 \end{equation}

Thus the stabilizer of an element $a_l\in\Fg(a_{l-1})$ equals
${\frak H}^l_a=\mathfrak{G}^l_a\cap\St$ in the case of the Lie
group, or
 $$
{\frak h}^l_a=\mathfrak{g}^l_a\cap\St
 $$
for the Lie algebra, where
 $$
\St=(\op{Ann}\tg_a)\circ_{\op{sym}}S^{l-1}T_a^*M\ot
T_aM+S^lT^*_aM\ot\tg_a.
 $$
In particular, since the symbol of 
$D^l(M)$ acts transitively, we get:
 \begin{equation}\label{zxc}
S^lT^*_aM\ot T_aM/\St\simeq S^l\tg^*_a\ot\vg_a.
 \end{equation}

 \begin{remark}\po\label{rk1}
The group $\mathfrak{G}^l_a$ for $l>1$ is abelian, which reflect
the affine property of the action, and so we can work only with
Lie algebras. In the case of 1-jets one should operate with the
Lie groups.
 \end{remark}

Now we specify our equivalence problem by a $G$-invariant
differential equation $\mathfrak{N}\subset J^l_r(M)$ on
submanifolds $N\subset M$ of codimension $r$. The symbol of this
equation $h_a^l\subset S^l\tg^*_a\ot\vg_a$ is a
$\rho_{l,l-1}$-vertical subspace of $T_{a_l}\mathfrak{N}$. Since
the pseudogroup $G$ acts on $\mathfrak{N}$, we obtain the
following exact sequence:
 \begin{equation}\label{130300}
0\to{\frak h}^l_a\hookrightarrow
\g^l_a\stackrel{\ll}\longrightarrow
h_a^l\stackrel{\varpi}\longrightarrow {\frak O}^l_a\to 0.
 \end{equation}

 \begin{dfn}\po
The quotient ${\frak O}^l_a=h_a^l/\ll(\g^l_a)$ is called the space
of \ $l$-{\em covariants\/} of the pseudogroup $G$ action. The
dual $({\frak O}^l_a)^*$ is named the space of \ $l$-{\em
variants\/}.
 \end{dfn}

Our study of formal equivalence of submanifolds under the
$G$-action is inductive and based on the following observation:

 \begin{prop}\po
Let $[N_1]^{l-1}_a=[N_2]^{l-1}_a\in\rho_{l,l-1}(\mathfrak{N})$ and
$l>1$. The $l$-jets of submanifolds $N_1$ and $N_2$ from
$\mathfrak{N}$ at a point $a\in M$ are $G$-equivalent if and only
if they belong to the same $\g^l_a$-orbit on $h_a^l$, which are
are affine subspaces of codimension equal $\dim{\frak O}_a^l=
\op{dim}h_a^l- \op{dim}\bigl(\g^l_a/{\frak h}^l_a\bigr)$. In other
words, this happens iff they have the same $l$-variants:
$\varpi([N_2]_a^l-[N_1]_a^l)=0$. \qed
 \end{prop}

The requirement $l>1$ is related to remark \ref{rk1}. For $l=1$
there is difference between symbolic Lie groups and algebras: In
the first case one gets orbits in the Grassmannian
$J^1_r(M)=\op{Gr}_r(T_aM)$, while in the latter one gets affine
subspaces in its tangent space at $a_1$. Thus 1-jets require a
separate treatment.

\subsection{\hpss Differential invariants and Tresse derivatives}\label{S13}

 \abz
Let $\mathcal{I}_k$ be the algebra of order $k$ differential
invariants of the pseudogroup $G$ action on $\mathfrak{N}_k$ (the
equation consists of pieces of different orders, see Appendix
\ref{A11}), i.e.\ functions constant on the $G_k$-orbits in
$\mathfrak{N}_k$. Denote by $\mathcal{I}$ the algebra of all
differential invariants. It is filtered by the subalgebras
$\mathcal{I}_k$ via the natural inclusion
$\rho_{k+1,k}^*:\mathcal{I}_k\to\mathcal{I}_{k+1}$ if the
pseudogroup $G$ is integrable. If the pseudogroup is not
integrable, we can still consider its finite piece to order $l$.

Sophus Lie proposed to produce new differential invariants via
invariant differentiations $\nabla$. He suggested a theorem that a
finite number of them $\nabla_1,\dots,\nabla_n$ is enough to
produce the whole algebra $\mathcal{I}$ from some $\mathcal{I}_k$.

An important case of invariant differentiations
$\nabla_i:\mathcal{I}_k\to\mathcal{I}_{k+1}$ constitute
derivatives a la Tresse, which we now introduce.

Suppose we have $n=\dim N=\dim M-r$ differential invariants
$f_1,\dots,f_n$ on ${\frak N}_k$. Provided $\pi_{k+1,k}({\frak
N}_{k+1})={\frak N}_k$ we define the differential operator
 $$
\hat\p_i:C^\infty({\frak N}_k)\to\Omega^1({\frak N}_{k+1}'),
 $$
where ${\frak N}_{k+1}'$ is the open set of points
$a_{k+1}\in{\frak N}_{k+1}$ with
 \begin{equation}\label{ndgTr}
df_1\we\dots\we df_n|_{L(a_{k+1})}\ne0.
 \end{equation}

We require that $\{f_i\}_{i=1}^n$ are such that ${\frak N}_{k+1}'$
is dense in ${\frak N}_{k+1}$. For the trivial equation ${\frak
N}_{k+1}=J^{k+1}_r(M)$ this is always the case. But if the
equation ${\frak N}$ is proper, this is a requirement of "general
position" for it. Given condition (\ref{ndgTr}) we write:
 $$
df|_{L(a_{k+1})}=\sum_{i=1}^n\hat\p_i(f)(a_{k+1})\,df_i|_{L(a_{k+1})},
 $$
which defines the function $\hat\p_i(f)$ uniquely at all the
points $a_{k+1}\in {\frak N}_{k+1}'$. This yields an invariant
differentiation $\hat\p_i=\hat\p/\hat\p
f_i:\mathcal{I}_k\to\mathcal{I}_{k+1}$. The expressions
$\hat\p_i(f)=\hat\p f/\hat\p f_i$ are called Tresse derivatives of
$f$ with respect to $f_i$.

The above construction can be presented more effectively in a
local chart $J^k\pi\subset J^k_r(M)$ (for this and the following
notions we refer to Appendix \ref{A1}). Given a local submanifold
$N\subset M$ we can find a transversal foliation of its
neighborhood and locally identify it with a bundle $\pi$ over $N$.
Then we can define Tresse derivative via the horizontal
differential $\hat d:C^\infty(J^k\pi)\to\Omega^1(J^{k+1}\pi)$.

In coordinate language given 1-jet $a_1=[N]_a^1$ we choose local
coordinates $(x^i,u^j)$ on $M$, with $\p_{x^i}$ being tangent to
$N$ at $a$ and $\p_{u^j}$ being transversal. Then $\hat d f=\sum
\D_i(f)dx^i$, where $\D_i$ is the operator of total derivative
with respect to coordinate $x^i$.

In these terms condition (\ref{ndgTr}) re-writes as:
  $$
\hat df_1\we\dots\we\hat df_n\ne0.
 $$
(in an open set $U$ -- a phrase we'll be omitting later on), i.e.
the Jacobian $\|\D_i(f_j)\|$ is non-degenerate. Then for any other
$f\in\mathcal{I}$ we have:
 \begin{equation}\label{hcoe}
\hat df=\sum_i\hat\p_i(f)\,\hat df_i.
 \end{equation}

Thus
 $$
\hat d=\sum dx^i\ot\D_{x^i}=\sum\hat df_i\ot\hat\p/\hat\p f_i,
 $$
which yields the expression of Tresse derivatives:
 \begin{equation}\label{dTresse}
\hat\p_i\stackrel{\text{def}}=\hat\p/\hat\p f_i=
\sum_j\Bigl(\D_{x^a}(f_b)\Bigr)^{-1}_{ij}\D_{x^j},
 \end{equation}
where $\bigl(\D_{x^a}(f_b)\bigr)$ is the Jacobian matrix in total
derivatives. This formula can be interpreted as a "change of
variables".

Informally speaking, $f_i$ are considered as base (horizontal)
coordinates on the equation $\mathfrak{N}$. They are classically
called differential parameters and in terms of them $\hat\p_i$ are
total derivatives. Then formula (\ref{hcoe}) has the standard
sense.

This idea was realized by S. Lie for vertical actions. This means
that the pseudogroup $G$ is represented in the equation ${\frak
N}\subset J^l\pi$ in such a way that every orbits in ${\frak N}_k$
belongs to a $\pi_k$-fiber. The base functions $x^1,\dots,x^n$
(for instance, local coordinates) are differential invariants. The
corresponding Tresse derivative $\hat\p_i$ coincides with the
operator of total derivative $\D_i$ with respect to coordinate
$x^i$.

Lie and his students believed this can fully extend to the general
pseudogroup actions and Tresse seems to be the first who realized
this.

\subsection{\hpss Covariants and equivalence}\label{S14}

 \abz
We will present now an infinitesimal analog of the construction of
differential invariants. Fix a point $a_l\in J^l_r(M)$ and define
the increasing filtration of $T_{a_l}^*J^l_r(M)$ by
 $$
\T_k(a_l)=\{d_{a_l}f\,|\,f\in\mathcal{I}_k\}\subset
T_{a_l}^*J^l_r(M),\quad k=0,\dots,l.
 $$

Note that $\T_l$ is the 1st order equation defining
$G^l$-differential invariants on $J^l_r(M)$ at regular points.
Near singular orbits the differential invariants have bad
behavior, and there we define the filtration as follows (the
definitions at regular points coincide):
 $$
\T_k(a_l)=\pi_{l,k}^*\op{Ann}T_{a_k}(G^k\cdot a_k).
 $$

 \begin{prop}\po
For $0<k\le l$: $\mathfrak{O}^k_a=(\T_k/\T_{k-1})^*$.
 \end{prop}

 \begin{proof}
In fact,
$\mathfrak{O}^k_a=T_{a_k}(\pi_{k,k-1})^{-1}_*(G^{k-1}\cdot
a_{k-1})/T_{a_k}(G^k\cdot a_k)\,$ and the claim follows.
 \end{proof}
 \vspace{4pt}

Since $(\mathfrak{O}^k_a)^*\subset S^k\tg_a\ot\vg^*_a$, we have
the natural map
 \begin{equation}\label{tresse}
\delta^*:(\mathfrak{O}^k_a)^*\ot\tg\to(\mathfrak{O}^{k+1}_a)^*,
 \end{equation}
which can be viewed as the symbol of invariant differentiation at
regular points. In order to prove surjectivity of this map for
large $k$, we will investigate the dual map and prove its eventual
injectivity, see \S\ref{S23}.

If we have $n$ independent differential invariants of order $k$,
then $\dim\mathfrak{O}^k\ge n$. In this case we can treat map
(\ref{tresse}) as an infinitesimal version of Tresse derivative.
This will provide a finite set of generators for differential
invariants, \S\ref{S24}.

Thus we get a solution to the formal equivalence problem by the
following inductive procedure. We start with a pseudogroup $G$ and
$\mathfrak{N}=J^l_r(M)$. Let the first nontrivial space of
$l$-covariants be $\mathfrak{O}^l$. Fix $l$-variants from
$(\mathfrak{O}^l_a)^*=\T^l$, i.e. fix order $l$ differential
invariants. If they are compatible as differential operators
(otherwise we need to add compatibility conditions), this yields a
smaller equation ${\mathfrak N}\subset J^l_r(M)$ on submanifolds
$N$ and we continue (in fact, the procedure is more complicated:
If the invariants are not constants, we take some of them as
"coordinates", express the others via them and fix the
corresponding functions-relations). At regular points the
procedure stops in a finite number of steps by the
Cartan-Kuranishi prolongation theorem.

An important case is an eventual absence of $l$-variants.

 \begin{dfn}\po\label{dfn5}
{\rm(i)} A pseudogroup $G$ is said to act {\em $l$-transitively\/}
near $a_l\in\mathfrak{N}$, if for any other jet
$b_l\in\mathfrak{N}$, close to $a_l$, there exists an element
$\vp_l\in G^l_{a,b}$ such that $\vp_l(a_l)=b_l$. In other words,
the orbit $G^l\cdot a_l$ is open. \\
{\rm(ii)} An action of a pseudogroup $G$ is said to be {\em
$l$-transversal\/} near $a_l$, if the above holds whenever
$a_{l-1}=b_{l-1}$. In other words, $\mathfrak{G}^l_a$ acts
transitively on $\mathfrak{F}(a_{l-1})$.
 \end{dfn}

To explain the word "transversality", consider the map
$\ll:\te\mapsto\bar\te$ from (\ref{lambda}). The space
$\ll^{-1}(h^l_a)\subset S^lT^*_aM\ot T_aM$ contains two subspaces
$\St$ and $\mathfrak{g}^l_a$.

Let $l>1$. The following statement follows from (\ref{zxc}),
(\ref{130300}) and definitions:

 \begin{prop}\po\label{thm4}
$l$-transversality of $G$ on $\mathfrak{N}$ is equivalent to any
of the conditions:
 \begin{itemize}
\item
 $\St$ is transversal to $\mathfrak{g}_a^l$ in $\ll^{-1}(h_a^l):$\
\ $\St+\mathfrak{g}_a^l=\ll^{-1}(h_a^l)$.
\item
 There are no $l$-covariants:
 $\mathfrak{O}_a^l=0$.\qed
 \end{itemize}
 \end{prop}

$l$-transversality is an inductive step to get $l$-transitivity.
Namely, we have:

 \begin{theorem}\po
Let $G^1\cdot a_1$ be open and $G$ acts $j$-transversally near
$a_j$ for $1<j\le l$. Then $G$ acts $l$-transitively near $a_l$.
\qed
 \end{theorem}

 \begin{dfn}\po
We will call an action of $G$ formally transitive if it is
$l$-transitive near a generic point of $\mathfrak{N}$ for every
$l$. If it is $l$-transversal for all $l$ starting from some
$l_0$, we will call such an action eventually transitive. This
basically means that the number of differential invariants is
finite.
 \end{dfn}

\section{\hps Homological methods}\label{S2}

 \abz
In this section we develop a technique to formally handle
differential invariants and prove the finiteness theorem.

\subsection{\hpss Cohomology of covariants}\label{S21}

 \abz
Consider a pseudogroup $G$ of order $k$. Denote
 $$
{\frak h}^{l-s,s}_a=\g^{l-s}_a\ot(\op{Ann}\tg_a\we
\La^{s-1}T^*_aM)+{\frak h}^{l-s}_a\ot\La^sT^*_aM
 $$
and let $\varrho:T^*_aM\to\tg^*_a$ be the restriction map.
Consider the following commutative diagram, where the horizontal
arrows are induced $\d$-differentials and the vertical ones are
obvious from exact four-sequence (\ref{130300}).
 $$
  \begin{CD}
 @. 0 @. 0 @. 0  \\
@. @VVV  @VVV @VVV @.\\
 0\!\! @>>> {\frak h}^l_a @>{\d}>>
 {\frak h}^{l-1,1}_a @>{\d}>>
 {\frak h}^{l-2,2}_a @>{\d}>>
 \!\!\dots \\
@. @VVV  @VVV @VVV @.\\
 0\!\! @>>> \g^l_a @>{\d}>>
 \g^{l-1}_a\ot T^*_aM @>{\d}>>
 \g^{l-2}_a\ot \La^2T^*_aM @>{\d}>>
 \!\!\dots \\
@. @V{\ll}VV  @V{\ll\ot\varrho}VV @V{\ll\ot\we^2\varrho}VV @.\\
 0\!\! @>>> \!\!h_a^l\!\! @>{\d}>>
 \!\!h_a^{l-1}\ot\tg^*_a\!\! @>{\d}>>
 \!\!h_a^{l-2}\ot\La^2\tg^*_a\!\! @>{\d}>>
 \!\!\dots \\
@. @VVV  @VVV @VVV @.\\
 0\!\! @>>> {\frak O}^l_a @>{\d}>>
 {\frak O}^{l-1}_a\ot\tg^*_a @>{\d}>>
 {\frak O}^{l-2}_a\ot\La^2\tg^*_a @>{\d}>>
 \!\!\dots \\
@. @VVV  @VVV @VVV @.\\
 @. 0 @. 0 @. 0  \\
  \end{CD}
 $$

Denote the cohomology of the first line at the term ${\frak
h}^{l-s,s}$ by $H^{l-s,s}({\frak h},\g)$ and the cohomology of the
forth line at the term ${\frak O}^{l-s}\ot\La^s\tg^*$ by
$H^{l-s,s}({\frak O})$. These latter will be called the cohomology
of covariants (in principle, they depend on the point of equation
$\mathfrak{N}$, but we will not indicate this).

The following statement is obtained by the usual diagram chase.

 \begin{prop}\po\label{lbl}
Suppose $H^{l-s-1,s+1}(\g)=H^{l-s-2,s+2}(\g)=0$ and
$H^{l-s,s}(h)=H^{l-s-1,s+1}(h)=0$. Then $H^{l-s,s}({\frak
O})\simeq H^{l-s-2,s+2}({\frak h},\g)$. \qed
 \end{prop}

 \begin{cor}\po
Let an order $k$ pseudogroup $G$ act on submanifolds $N\subset M$
of fixed codimension $r$, more precisely on $J^k_r(M)$. Let $G$ be
$(q+2)$-acyclic and $l>k+1$. Then $H^{l-s,s}({\frak O}) \simeq
H^{l-s-2,s+2}({\frak h},\g)$ for all $s\le\min(l-k-2,q)$. In
particular, if $G$ is involutive, then the equality holds for all
$s\le l-k-2$. \qed
 \end{cor}

 \begin{cor}\po
Consider a $(q+2)$-acyclic pseudogroup $G$ of order $k$ acting on
an equation ${\mathfrak N}\subset J^m_r(M)$, which is
$(p+1)$-acyclic. Let $l>\max(k+1,m)$. Then $H^{l-s,s}({\frak
O})\simeq H^{l-s-2,s+2}({\frak h},\g)$ for
$s\le\min(l-k-2,l-m-1,p,q)$. \qed
 \end{cor}

 \begin{cor}\po\label{cor21}
Suppose that: 1) ${\mathfrak O}_a^{l-1}=0$; 2)
$h_a^l=(h_a^{l-1})^{(1)}$, $\g_a^l=(\g_a^{l-1})^{(1)}$; 3)
$H^{l-2,2}(\g_a)=0$. Then ${\mathfrak
O}_a^l=H^{l-2,2}(\mathfrak{h}_a,\g_a)$. \qed
 \end{cor}

Thus we obtain a method to calculate recursively the space of
covariants $\mathfrak{O}^l_a$ if we know the cohomology groups
$H^{*,*}({\frak h},\g)$. This leads to the inductive approach of
\S\ref{S14} to the equivalence problem. Due to proposition
\ref{lbl}:
 $$
H^{l,0}({\frak h},\g)=0,\qquad H^{l,1}({\frak h},\g)=0 \text{ for
} l\ge k.
 $$
Let us calculate the groups $H^{l,s}({\frak h},\g)$ for $s>1$. We
do it at first with an additional assumption of
non-characteristisity.

 \begin{theorem}\po\label{th22}
Let $G$ be $q$-acyclic and let $c=\min(l-k,q)$. Denote by
$H^{l-s,s}({\frak h})$ the cohomology group of the complex
  \begin{equation}\label{hhh2}
 0 \to{\frak h}^l_a \to {\frak h}^{l-1}_a\ot\tg^*_a \to
{\frak h}^{l-2}_a\ot\La^2\tg^*_a \to \dots
  \end{equation}
at the term ${\frak h}^{l-s}_a\ot\La^s\tg^*_a$.

Suppose that the subspace $\tg_a\subset T_aM$ is strongly
non-characteristic for $\g_a$ \cite{KL3}, i.e.
$\op{Ann}(\tg_a)\circ S^{k-1}T^*_aM\ot T_aM\cap\g_a=0$. Then for
$0\le s<c$ we have:
 $$
H^{l,s}({\frak h},\g)=H^{l,s}({\frak h}).
 $$
 \end{theorem}

 \begin{proof}
Consider the following commutative diagram of vertical exact
three-sequences, where $\d'$ is the induced differential:
 $$
  \begin{CD}
 @. 0 @. 0 @. 0  \\
@. @VVV  @VVV @VVV @.\\
 0 @>>> 0 @>>>
 \!\!\g^{l-1}_a\ot\op{Ann}\tg_a\!\! @>{\d}>>
 \!\!\g^{l-2}_a\ot\op{Ann}\tg_a\we T^*_aM\!\! @>{\d}>>
 \dots \\
@. @VVV  @VVV @VVV @.\\
 0 @>>> \!\!\g^l_a\!\! @>{\d}>>
 \g^{l-1}_a\ot T^*_aM @>{\d}>>
 \g^{l-2}_a\ot\La^2T^*_aM @>{\d}>>
 \dots \\
@. @VVV  @VVV @VVV @.\\
 0 @>>> \!\!\g^l_a\!\! @>{\d'}>>
 \g^{l-1}_a\ot\tg^*_a @>{\d'}>>
 \g^{l-2}_a\ot\La^2\tg^*_a @>{\d'}>>
 \dots \\
@. @VVV  @VVV @VVV @.\\
 @. 0 @. 0 @. 0  \\
  \end{CD}
 $$
The middle line is $c$-acyclic. If $\tg_a$ is strongly
non-characteristic, we have the same property for the bottom line
\cite{KL3}. Let $H^{i,j}(\g\ot\op{Ann}\tg;\d)$ denote the
cohomology of the first complex at the term
$\g_a^i\ot\op{Ann}\tg_a\we\La^jT^*_aM$. A diagram chase gives:
$H^{i-1,j}(\g\ot\op{Ann}\tg;\d)\simeq H^{i,j}(\g,\d')=0$ for
$i>k$, $0\le j<c$.

Consider the following commutative diagram with vertical
three-sequences being exact. Note that if the subspace $\tg_a$ is
strongly non-characteristic, we can consider ${\frak h}^l_a\subset
S^l\tg^*_a\ot\tg_a$, so that the bottom complex is the usual
Spencer $\d$-complex on $\tg_a$.
 $$
  \begin{CD}
 @. 0 @. 0 @. 0  \\
@. @VVV  @VVV @VVV @.\\
 0 @>>> 0 @>>>
 \!\!\g^{l-1}_a\ot\op{Ann}\tg_a\!\! @>{\d}>>
 \!\!\g^{l-2}_a\ot\op{Ann}\tg_a\we T^*_aM\!\! @>{\d}>>
 \dots \\
@. @VVV  @VVV @VVV @.\\
 0 @>>> \!\!{\frak h}^l_a\!\! @>{\d}>>
 {\frak h}^{l-1,1}_a @>{\d}>>
 {\frak h}^{l-2,2}_a @>{\d}>>
 \dots \\
@. @VVV  @VVV @VVV @.\\
 0 @>>> \!\!{\frak h}^l_a\!\! @>{\d'}>>
 {\frak h}^{l-1}_a\ot\tg^*_a @>{\d'}>>
 {\frak h}^{l-2}_a\ot\La^2\tg^*_a @>{\d'}>>
 \dots \\
@. @VVV  @VVV @VVV @.\\
 @. 0 @. 0 @. 0  \\
  \end{CD}
 $$

Since the first horizontal complex is $c$-acyclic, the middle and
the bottom complexes have the same cohomology in the first $c$
terms.
 \end{proof}
\vspace{4pt}

In the Spencer complex on $\tg_a$ all $\d$-cohomology groups
eventually vanish (Poincar\'e $\d$-lemma \cite{S,KLV}). In
non-characteristic case for a big number $i$ (actually such big
that the equation $G$ on $TM$ and its restriction to $\tg$ as well
as the equation ${\mathfrak N}$ on $\tg$ become involutive) we
have: $H^{i,j}({\frak O})=0$ .

Thus if the pseudogroup does not have all subspaces of given
codimension $r$ weakly characteristic \cite{KL3}, then we have the
following finiteness theorem: Cohomology of covariants eventually
vanish (on an open dense subset of the equation ${\mathfrak N}$).
We will prove in \S\ref{S23} that this is a general fact.

Note however that with the approach of this section we calculated
the cohomology of covariants, which is an important invariant of
pseudogroup action:
 \begin{cor}\po
Suppose that assumptions of Proposition \ref{lbl} and Theorem
\ref{th22} hold. Then $H^{l-s,s}({\frak O})\simeq
H^{l-s-2,s+2}({\frak h})$ for $0\le s<c$. \qed
 \end{cor}

\subsection{\hpss Criterion of transversality}

 \abz
By Corollary \ref{cor21} a very important cohomology group of
${\mathfrak h}$ is $H^{l,2}({\frak h})$.

 \begin{theorem}\label{H-clas}\po
Let a pseudogroup $G$ be 2-acyclic: $H^{l,2}(G)=0$, $l\ge k$.
Suppose that for some number $l_0>k$ the submanifold $N$ at a
point $a$ is $l_0$-transversal with respect to the pseudogroup $G$
action. Assume also that $H^{l,2}({\frak h})=0$ and
$H^{l,1}(\g,\d')=0$ for $l>l_0$. Then $N$ is $l$-transversal for
all $l>l_0$ at $a$.
 \end{theorem}

 \begin{proof}
Indeed, from the first diagram of the proof of Theorem \ref{th22}
we get: $H^{l-1,1}(\g\ot\op{Ann}\tg;\d)\simeq H^{l,1}(\g,\d')=0$.

From the second diagram of the same proof we obtain that since
$H^{l-1,1}(\g\ot\op{Ann}\tg;\d)=0$ the map of cohomology
$H^{l-2,2}({\frak h};\g)\to H^{l-2,2}({\frak h})$ is injective.
Thus by assumptions and corollary \ref{cor21}: ${\frak
O}^l_a\simeq H^{l-2,2}({\frak h}_a;\g_a)=0$ for all $l\le l_0$.
 \end{proof}

 \begin{cor}\po
With the assumptions of Theorem \ref{H-clas} the pseudogroup
action is eventually transitive. If we assume in addition that $N$
is $l$-transversal with respect to the pseudogroup $G$ action for
all $l<l_0$ and that the orbit $G^1\cdot a_1$ is open, then the
action of $G$ is formally transitive around $a$. \qed
 \end{cor}
 \vspace{4pt}

Note that in the theorem we don't require $\tg_a$ to be strongly
non-characteristic. This means that zero cohomology of the bottom
complexes from diagrams in Theorem \ref{th22} can be non-vanishing
even for large $l$.

However often the other cohomology groups vanish in stable range
(big $l$). This is related to the following fact:

 \begin{prop}\po
Let the pseudogroup $G$ be 2-acyclic from some level $l_0$.
Suppose that $H^{l,1}(\g_a,\d')=0$ for $l\ge l_0$. Then
$H^{l,1}({\frak h}_a)=0$.
 \end{prop}

 \begin{proof}
From the first diagram of Theorem \ref{th22} we get the
isomorphism $H^{l-1,1}(\g\ot\op{Ann}\tg;\d)\simeq
H^{l,1}(\g;\d')=0$.

From the second diagram since $H^{l-1,1}(\g\ot\op{Ann}\tg;\d)=0$
we obtain that the map of cohomology $H^{l,1}({\frak h};\g)\to
H^{l,1}({\frak h})$ is surjective. The claim follows from the fact
that $H^{l,1}({\frak h};\g)=0$.
 \end{proof}
 \vspace{4pt}

This means that the complex (\ref{hhh2}) is natural in the
following sense:
 $$
{\frak h}^l=\{\te\in\g^l\,|\,\p_v(\te)\in{\frak h}^{l-1}\,\forall
v\in\tg_a\}.
 $$

Notice that there exists an important necessary condition for
eventual (and hence formal) transitivity of the pseudogroup
action. This is a purely dimensional obstruction to
transversality.

Namely, by proposition~\ref{thm4} $l$-transversality condition
imposes the following inequality on the symbol $h_a^l$ of the
equation $\mathfrak{N}$:
 \begin{equation}\label{dim-trans}
\op{dim}\g^l_a\ge\op{dim}h_a^l.
 \end{equation}
This easy-to-check condition is often helpful. Namely, in many
cases its fulfilment implies transversality for {\em generic\/}
submanifolds $N$ (see examples below).

\subsection{\hpss Finiteness theorem}\label{S23}

 \abz
Here we prove an algebraic point-wise version of the finiteness
theorem. Its local version will appear in the next section.

For a Lie pseudogroup $G$ the corresponding Lie equation for
vector fields is linear. Then the characteristic variety
$\op{Char}^\CC(G;\vp_l)$ (we refer to Appendix \ref{A1} for the
definition and properties) of $G^l$ depends only on the base point
$a=\rho_{l,0}(\vp_l)$.

More generally, the same holds for any pseudogroup $G$ after some
number of prolongations, i.e. for some $l\ge l_0$. Indeed, if
$G^l$ has prolongation over points
$\vp_l',\vp_l''\in\rho_{l,0}^{-1}(a)$ (this is given by the
conditions $W_j(G;\vp_j')=W_j(G;\vp_j'')=0$), then the
characteristic varieties
$\op{Char}^\CC(G;\vp_j'),\op{Char}^\CC(G;\vp_j'')\subset P^\CC
T^*_aM$ coincide. We will denote the characteristic variety also
by $\op{Char}^\CC(\g_a)$.

 \begin{theorem}\po\label{vanish}
The Grassmannian space $J^1_r(M)_a=\op{Gr}_n(T_aM)$ of
$n$-dimensional subspaces $\tg_a\subset T_aM$, $n+r=m=\dim M$,
contains an open dense subset $\uu(a)$, depending only on
$\op{Char}^\CC(\g_a)$, and there exists a number $l_0$, depending
only on the pseudogroup $G$ and the equation for submanifolds
${\frak N}$, such that the following holds. For any point
$a_l\in{\frak N}_l$, $l\ge l_0$, with $a=\rho_{l,0}(a_l)$ and such
that $a_1=\rho_{l,1}(a_l)$ is an admissible tangent space
$\tg_a\in\uu(a)$ we have:
 $$
H^{i,j}({\frak O})=0\quad \text{ for any }\ i+j=l\ge l_0.
 $$
 \end{theorem}

 \begin{proof}
The proof of the theorem is split into two two parts, depending on
weather $r\le\op{codim}\op{Char}^\CC(\g)$ or
$r\ge\op{codim}\op{Char}^\CC(\g)$ (in the case of equality both
approaches are equivalent).

Note that the codimension does not depend on weather we consider
affine variety in $T_x^\CC M$ or its projectivization in $PT_x^\CC
M$ ($P$ denotes projectivization and $^\CC$ -- complexification).
However usage of complex characteristics is crucial. Also note
that we do not require the characteristic variety to be
irreducible, but take $d=\op{codim}\op{Char}^\CC(\g)$ to be the
codimension of its regular component (so this value is the minimum
of codimensions by all regular points of all irreducible pieces).
We have: $d\in[0,n]$.

{\bf 1.} $r\le d$. In this case $\uu(a)$ consists of subspaces
$\tg$ such that $P\op{Ann}(\tg)^\CC$ does not intersect
$\op{Char}^\CC(\g)$. It is possible by Noether normalization lemma
and all generic subspaces $\tg$ are such.

For a vector space $V$ denote by $SV=\oplus S^iV$ the ring of
homogeneous polynomials on $V^*$. Let $I_0(\g)$ be the annihilator
of the subvariety $\op{Char}^\CC(\g)\subset P(T_a^*M)^\CC$, i.e.
the ideal of homogeneous polynomials vanishing on the
characteristic variety. It equals the radical of the
characteristic ideal $I(\g)\subset S(T_a^\CC M)$ (see
Appendix~\ref{A1}); here again $S(T_a^\CC M)=\oplus S^i(T_a^\CC
M)$ is the polynomial algebra.

In addition Noether lemma states \cite{M} that the projection
along annihilator $P_\tg:\op{Char}^\CC(\g)\to P(\tg^*)^\CC$ is a
finite-to-one closed map such that the homogeneous ring $S(T_a^\CC
M)/I_0(\g)$ is a finitely generated module over the algebra
$S(\tg^\CC)$.

We claim that the homogeneous ring $S(T_a^\CC M)/I(\g)$ is a
finitely generated module over the algebra $S(\tg^\CC)$. Indeed,
let us have a polynomial relation in the ring $S(T_a^\CC
M)/I_0(\g)$:
 $$
Q(f_1,\dots,f_m)\in I_0(\g),\quad f_i\in S(T_a^\CC M),\ Q\in
S(\tg^\CC).
 $$
Denote by $N$ the minimal integer number such that
$I_0(\g)^N\subset I(\g)$. Then we have the following polynomial
relation in the ring $S(T_a^\CC M)/I(\g)$:
 $$
Q^N(f_1,\dots,f_m)\in I(\g).
 $$

Thus the characteristic module $\g^*$, dual to the symbolic system
$\g$ (see Appendix \ref{A1}), is Noetherian over $S(\tg)$
(informally: the symbolic module grows over the characteristic
variety and it is projected finite-to-one). Thus the Koszul
cohomology of $\g^*$ is finite. Dualization yields finiteness of
the Spencer cohomology $H^{*,*}(\g,\d')$ of $\g$ over $\tg$.

Alternatively the latter claim follows from Poincar'e $\d$-lemma
\cite{S,KLV}. The bound $l_0$ such that $H^{i,j}(\g,\d')=0$ for
$i\ge l_0$ depends only on dimensions of the module $\g^*$ over
the algebra $S(\tg)$ and so is universal over all $\tg\in\uu(a)$.

Let us take $l_0$ such that the $l^\text{th}$ Spencer
$\d$-complexes for $\g$ and $h$ are acyclic, when $l\ge l_0$. For
such $l$ in the first commutative diagram from the proof of
Theorem \ref{th22} the second and the third complexes are acyclic.
Therefore the first one is acyclic. It is also the first complex
of the second diagram from the proof, so that we get isomorphism
between the cohomology of the second and third complexes, i.e.
 $$
H^{i,j}({\frak h},\g)\simeq H^{i,j}({\frak h})\quad\text{ for }\
i+j=l\ge l_0.
 $$

From the above isomorphism we deduce vanishing of the zero and
first cohomology of complex (\ref{hhh2}) in the range $l\ge l_0$.
This means that for these $l$ the space ${\frak h}^{l+1}$ is the
Spencer $\delta$-prolongation of the space ${\frak h}^l$. This
implies (again by Poincar\'e $\d$-lemma) that the cohomology
$H^{i,j}({\frak h})$ vanish for big $i+j=l$. Thus increasing $l_0$
properly, we obtain that the cohomology $H^{i,j}({\frak h},\g)=0$
for $i+j=l\ge l_0$.

Now the claim follows from Proposition \ref{lbl} because the
cohomology of covariants coincide with the cohomology
$H^{*,*}({\frak h},\g)$ in the stable range $l\ge l_0$.

{\bf 2.} $r\ge d$. In this case $\uu(a)$ consists of subspaces
$\tg$ such that the projection of $\op{Char}^\CC(\g)$ along
$P\op{Ann}(\tg)^\CC$ on $P(\tg^*)^\CC$ is surjective. Again all
generic subspaces $\tg$ are such due to Noether normalization
lemma \cite{M}.

An element $v\in T_aM$ is regular (in the sense of commutative
algebra \cite{AB,BH}) if it does not belong to the annihilator
$I(\g)$ of the module $\g^*$. This means that $P\op{Ann}(v)^\CC$
does not contain the characteristic variety $\op{Char}^\CC(\g)$.
This is equivalent to the fact that the projection of
$\op{Char}^\CC(\g)$ along $P\op{Ann}(v)^\CC$ to $P(\CC v)^*$ is
not empty and is therefore surjective.

More generally, a sequence $(v_1,\dots,v_n)$ is regular
($\g^*$-sequence) iff the projection of $\op{Char}^\CC(\g)$ along
$P\op{Ann}(v_1,\dots,v_n)^\CC$ to $P(\langle
v_1,\dots,v_n\rangle^\CC)^*$ is surjective. We conclude that there
exists a regular sequence $(v_1,\dots,v_n)$ in $\tg$ of length
$n=m-r=\dim\tg\le m-d$.

This implies that all the Koszul homology of the module $\g^*$
w.r.t. the sequence $(v_1,\dots,v_n)$, or equivalently with
coefficients in $\tg$, vanish except for the zero cohomology
group, see \cite{AB} or the appendix (including a letter of Serre)
in \cite{GS1} (equivalently we can say that $\g^*$ is a
Cohen-Macaulay module over $S(\tg)$, which implies the same result
\cite{AB,BH}). Dualizing this statement we obtain that the Spencer
cohomology groups vanish: $H^{i,j}(\g,\d')=0$, $i\ge k$, $0<j\le
n$. The zero cohomology group $H^{i,0}(\g,\d')$ for $r>d$ is
always non-zero and can be non-zero even for $r=d$.

As in the first case we use two diagrams from Theorem \ref{th22}
to conclude that the second and the third complexes of the second
diagram have the same cohomology, save for the zero cohomology
(which is zero for the second complex, but can be non-zero for the
third one):
 $$
H^{i,j}({\frak h},\g)=H^{i,j}({\frak h})\quad\text{ for }\ i\ge
k,\ j>0.
 $$
In particular, $H^{i,1}({\frak h})=0$ for $i\ge k$. This again
yields that the positive cohomology of ${\frak h}$ eventually
vanish: $H^{i,j}({\frak h})=0$ for $i+j=l\ge l_0$, $j>0$.
Consequently, $H^{i,j}({\frak h},\g)=0$ for $i+j=l\ge l_0$ and all
$j$.

Applying Proposition \ref{lbl} we again get vanishing of the
cohomology of covariants $H^{i,j}({\frak O})$ in the stable range
$i+j=l\ge l_0$.
 \end{proof}
\vspace{4pt}

Notice that with the approach of Theorem \ref{vanish} the estimate
for the place, where cohomology vanish, can be much higher than
that one of Theorem \ref{th22}. However the latter case works only
for pseudogroups such that not all subspaces are weakly
characteristic. But for some important pseudogroups, like
volume-preserving or symplectic pseudogroups, all tangent
subspaces are weakly characteristic. The finiteness theorem
however still holds even in such cases.

 \begin{remark}\po
There is another approach to prove Theorem \ref{vanish}. Namely,
for big $l_0$ the symbolic system $\{\g^l\}_{l\ge l_0}$ is
involutive, so that all Spencer $\d$-cohomology groups vanish.
This means that $S(T_aM)$-module $\oplus_{l\ge l_0}(\g^l)^*$
is Cohen-Macaulay. Then almost every subspace $\tg$ contains a
regular sequence and then its positive Koszul homology vanish, so
that $H^{i,j}(\g,\d')=0$ for $i\ge l_0$ and $j>0$.

With this approach we however cannot explicitly formulate which
subspaces $\tg$ are good for vanishing of positive cohomology of
the complex
 $$
0\to\g^l_a\stackrel{\d'}\longrightarrow\g^{l-1}_a\ot\tg^*_a
\stackrel{\d'}\longrightarrow\g^{l-2}_a\ot\La^2\tg^*_a\to\cdots.
 $$
This is the crucial place in the proof and the rest is just the
diagram chase.
 \end{remark}

\subsection{\hpss Relation to the theorems of Lie and
Tresse}\label{S24}

 \abz
Let us formulate now the regularity assumptions. We let the point
$a_l$ vary over ${\frak N}_l$ with big enough $l\ge l_0$, so that
the ranks of the symbol bundles are locally constant,
$\d$-cohomology are stabilized etc.

We call a point $a_l$ regular if the space
$\tg_a=a_1=\rho_{l,1}(a_l)$ is admissible in the sense of Theorem
\ref{vanish}. The collection of such points is open and will be
denoted by:
 $$
\op{Reg}_l^1({\frak N},G)=\{a_l\in{\frak N}_l\,|\, a_1\in\uu(a)\}.
 $$
We want to claim that $\op{Reg}^1_{l_0}({\frak N},G)$ is dense in
${\frak N}_{l_0}$. This is so if the equation ${\frak N}$ is
trivial -- defined by empty set of relations, i.e. ${\frak
N}_l=J^l_r(M)$.

More generally, each equation ${\frak N}$ with sufficiently rich
${\frak N}_1$ is such, meaning that any jet $a_1$ can be perturbed
to $a_1'$ in the fiber over $a=\rho_{1,0}(a_1)$ to satisfy the
transversality conditions of Theorem\ref{vanish}: $a_1'\in\uu(a)$.
In other words, ${\frak N}_1$ is not contained in the
$G$-invariant singular equation $\cup_{a\in M}
[\op{Gr}_n(T_aM)\setminus\uu(a)]$.

We will need another assumption, which is similar to Kumpera's
hypothesis H$_3$ \cite{Kum}. Denote by $\Delta_l(a_l)$ the tangent
space to the $G_l$-orbit through $a_l\in{\frak N}_l$. Recall that
each $a_{l+1}\in{\frak N}_l^{(1)}$ determines a horizontal space
$L(a_{l+1})\subset T_{a_l}{\frak N}_l$. Consider the open set
 $$
\op{Reg}_l^2({\frak N},G)=\{a_l\in{\frak N}_l\,|\, \exists
a_{l+1}\in{\frak N}_{l+1}: \Delta_l(a_l)\cap L(a_{l+1})=0\}.
 $$
We want to claim that $\op{Reg}_{l_0}^2({\frak N},G)$ is also
dense. This can fail, for instance, if we have $l$-transversality.
For the trivial equation ${\frak N}_l=J^l_r(M)$ the condition
means that we have $n$ independent differential invariants
$f_1,\dots,f_n$. However for proper ${\frak N}$ this is a
requirement on the equation.

Given these two regularities we will prove that the algebra of
differential invariants has a finite base w.r.t. Tresse
derivatives (on an open dense set; but this condition is natural
since usually the differential invariants have singularities).

S. Lie used invariant differentiations to generate differential
invariants. A. Tresse observed that they can be obtained if we
have a sufficient number of independent differential invariants.

 \begin{remark}\po\label{remk3}
A. Kumpera proved Lie-Tresse theorem for a Lie sheaf of vector
fields \cite{Kum}. Under his conditions the maps
$\lambda\ot\La\varrho$ in the four-line diagram of \S\ref{S21}
are injective, so that the first complex of it vanishes. Then the
diagram becomes with 3-lines and after regularity assumptions the
vanishing theorem follows from the stabilization of the cohomology
of pseudogroup $G$ as Cartan-Kuranishi theorem \cite{Kur} states.
 \end{remark}

We will deduce now the theorem of Lie-Tresse. We will assume at
first that both equations $G$ of the pseudogroup and ${\frak N}$
for the submanifolds are formally integrable.

 \begin{theorem}\po
Let a pseudogroup $G$ act on an equation ${\frak N}$. Suppose that
both are formally integrable. Let also $\op{Reg}_{l_0}^1({\frak
N},G)$ and $\op{Reg}_{l_0}^2({\frak N},G)$ be dense in ${\frak
N}_{l_0}$. Then the infinitely prolonged equation ${\frak
N}^{(\infty)}$ contains a (no more than countable) collection of
open $G$-invariant sets $U_\alpha$, the union of which $U=\cup_\a
U_\a$ is dense, with the following properties.

Consider some $U_\a$. Then there are $n$ differential invariants
$f_1,\dots,f_n$ on it, with the corresponding invariant
differentiations $\hat\p_1,\dots,\hat\p_n$, and some other
differential invariants $g_1,\dots,g_m$ such that all differential
invariants in $U_\a$ can be expressed via the $g_j$ and their
invariant derivatives $\hat\p^J(g_j)$ (for a multi-index
$J=(j_1,\dots,j_n)$ we denote
$\hat\p^J=\hat\p_1^{j_1}\cdots\hat\p_n^{j_n}$).
 \end{theorem}

Usually (so-called regularity assumptions) there is only one such
set $U$.
 \vspace{4pt}

 \begin{proof}
By the assumption $U(l_0)=\op{Reg}_{l_0}^1({\frak N},G)\cap
\op{Reg}_{l_0}^2({\frak N},G)$ is dense in ${\frak N}_{l_0}$. If
${\frak N}_l={\frak N}_{l_0}^{(l-l_0)}$, then the characteristic
variety on the level $l$ is the same as on the level $l_0$. Thus
$\op{Reg}^1_l=\rho_{l,l_0}^{-1}(\op{Reg}^1_{l_0})\cap{\frak N}_l$.
The same applies to $\op{Reg}^2_l$. Thus we let
$U=\rho_{\infty,l_0}^{-1}(U(l_0))$ and this set can be represented
as a union of sets $U_\a=\rho_{\infty,l_0}^{-1}(U_\a(l_0))$.

Without loss of generality we suppose that each $U_\a$ belongs to
a local chart, so that we can work with the jets of a bundle
$J^l\pi$, which has a convenient representation (\ref{dTresse}).
Moreover since $U_\a(l)\subset \op{Reg}_l^2({\frak N},G)$ for
$l\ge l_0$ there are $n=\dim\tg$ differential invariants
$f_1,\dots,f_n$ such that $\hat df_1\we\dots\we\hat df_n\ne0$ (we
can shrink $U_\a$). Thus the Tresse derivatives
$\hat\p_i=\hat\p/\hat\p f_i$ are well-defined in $U_\a$.

Let us calculate the symbols of these differentiations. For this
we need a lift of vectors to invariant differentiations, described
below.

Consider a jet $a_{l+1}\in{\frak N}_{l+1}$,
$a_l=\pi_{l+1,l}(a_{l+1})$, $l\ge l_0$. Restrictions of $\hat
df_i$ to the horizontal plane $L(a_{l+1})\subset T_{a_l}{\frak
N}_l$ form a basis. Thus we have a basis $e_1^*,\dots,e_n^*$ of
$\tg_a^*$ given by $\pi_l^*e_i^*=\hat df_i|_{L(a_{l+1})}$.

Denote by $e_1,\dots,e_n$ the dual basis of $\tg_a$ and let
$v\in\tg_a$. Decompose $v=\sum_{i=1}^nv_ie_i$. Choose a system of
local coordinates $(x^i)_{i=1}^n$ near $a\in M$ such that
$e_i=\p_{x^i}$ at $a$. Then the symbol map
 \begin{equation}\label{311}
({\frak O^l_a})^*\ot\tg_a\longrightarrow({\frak O^{l+1}_a})^*
 \end{equation}
is given by the formula
 $$
[d_{a_l}f]\ot v\mapsto
\Bigl[\sum\nolimits_{i=1}^nv_i\D_{x^i}(f)\Bigr].
 $$
Here $[d_{a_l}f]$ represents an element of $({\frak
O^l_a})^*\subset S^l\tg_a\ot\vg_a^*$ and in the right-hand-side we
restrict the covector at the point $a_{l+1}$ to the vertical
subspace $T_{a_{l+1}}^\text{vert}{\frak N}_{l+1}\subset
S^{l+1}\tg_a^*\ot\vg_a$ and take the quotient. The result does not
depend on the coordinate system $(x^i)_{i=1}^n$ adapted at $a$ as
indicated above.

Recall that $({\frak O}^{l+1}_a)^*=
\{df\,|\,f\in\mathcal{I}_{l+1}\}/ \{df\,|\,f\in\mathcal{I}_l\}$.
Since $l\ge l_0$ and $U_\a(l)\subset \op{Reg}_l^1({\frak N},G)$,
the cohomology of covariants vanish, which means that map
(\ref{311}) is epimorphic. This implies that
 $$
\mathcal{I}_{l+1}=\langle\rho^*_{l+1,l}(\mathcal{I}_l),
\hat\p_1(\mathcal{I}_l),\dots,\hat\p_n(\mathcal{I}_l)\rangle.
 $$
Indeed, by the finiteness theorem the differentials of the
functions on the right span the whole space of differentials of
the functions to the left. The claim follows from the implicit
function theorem.
 \end{proof}

 \begin{remark}\po
We explained the density condition of $\op{Reg}_{l_0}^1({\frak
N},G)$ before Remark \ref{remk3}. In certain cases density of
$\op{Reg}_{l_0}^2({\frak N},G)$ in ${\frak N}_{l_0}$ can be also
guaranteed.

Indeed, this is so if the dimension of characteristic variety of
the equation ${\frak N}$ exceeds the dimension of characteristic
variety of the pseudogroup $G$ (just by comparison of Poincar\'e
polynomials for the corresponding symbolic modules). This latter
condition is realized for infinite pseudogroups $G$ acting on
jet-spaces $J^\infty(\pi)$ with "functional dimension"
$\op{dim}_\CC\op{Char}^\CC(G)<\op{rank}(\pi)$ and for
(finite-dimensional) Lie groups $G$ acting on equations ${\frak
N}$ of infinite type.
 \end{remark}

We will now argue that the theorem holds for general
non-integrable (for both $G$ and ${\frak N}$) case as well and
show how the equation ${\frak N}$ on submanifolds naturally
appears. This is contained in the following three remarks.

{\bf 1.} When we consider the orbits of the pseudogroup $G$ even
in the pure jet-space $J^l_r(M)$ (but maybe on equation) there are
regular and singular orbits. The setup for constructing invariant
differentiations requires to restrict to the former. In addition,
the differential invariants as well as the Tresse derivatives are
usually  not defined on the whole space (for instance, because
these absolute differential invariants can be obtained as ratio of
two relative differential invariants). Thus we need to remove a
closed nowhere dense subset of the jet-space. This subset (or its
regular part) is a $G$-invariant equation ${\frak N}$. And we need
to apply the machinery to it. In turn, the orbits in it are
divided into regular and singular, so that we get smaller equation
etc.

{\bf 2.} If the equation ${\frak N}$ or the pseudogroup $G$ are
not integrable, then we must use the prolongation-projection
scheme. Each time we obtain a set of compatibility conditions we
project the equation $G_k$ (or ${\frak N}_k$) to obtain new
equations of smaller order, which we prolong etc. When this
concerns the pseudogroup, the space of differential invariants
grows: To the existent invariants we add new. In addition, we can
preserve existing invariant differentiations if we already possess
some, adding only new generators -- differential invariants
$g_m,\dots,g_{m+s}$.

Note that a shrink of $G$ results in a shrink of the
characteristic variety $\op{Char}^\CC(G)$ and a shrink of ${\frak
N}$ results in a shrink of its prolongation. Therefore both
regularity sets $\op{Reg}_l^1$ and $\op{Reg}_l^2$ can change and
we should care that the density property for regular points is not
lost. However by Cartan-Kuranishi theorem there is only a finite
number of such shrinks in the process of prolongation-projection.
On each of this step we add a finite number of differential
invariants, which remain invariants during the rest of the
process.

{\bf 3.} The Tresse derivatives as introduced in \S\ref{S13} are
invariant differentiations provided that $G$ and ${\frak N}$ are
integrable (this is not obvious from (\ref{dTresse}), but follows
from the preceding formulas). But since we arrive to integrable
equations in a finite number of steps, we will eventually get the
required differentiations or observe finiteness of invariants.
Thus even in non-integrable (but sufficiently regular) case Lie
theorem holds.

\section{\hps Invariants of geometric structures}

 \abz
In this section we check the transitivity condition for the
automorphisms pseudogroups of some basic geometric structures.
Irreducible Lie pseudogroups were classified by E.\,Cartan
(\cite{C1}). We consider at first these integrable pseudogroups.
Examples of this section were mostly considered in our preceding
paper \cite{KL2}, so the proofs will be omitted, though we
indicate how to obtain the results from our cohomological
machinery.

Note that the $l$-pseudogroup $G^l$ consists of the jets of
diffeomorphisms preserving the structure to order $l$. So if the
structure is non-integrable, then the prolongation-projection
method changes the equation and the sub-pseudogroup $G^j$ can be
different as embedded into $G^k$ and $G^l$, $j<\min(k,l)$ (so one
should be careful with notations). We consider examples of the
transformation pseudogroups of non-integrable structures at the
end of the section.

In this section we suppose ${\frak N}=J^l_r(M)$ unless the
contrary is stated (thus $h_a^l=S^l\tg_a^*\ot\vg_a$).

\subsection{\hpss General and volume preserving pseudogroups}\label{GLSL}

 \abz
The general pseudogroup $G=\op{Diff}_\text{loc}(M)$ is involutive.
We have: $G^l=D^l(M)$ for all $l$. In this case ${\frak
h}^l=(\op{Ann}\tg)\circ_{\text{sym}}S^{l-1}T^*M\ot TM+
S^lT^*M\ot\tg$. Complex (\ref{hhh2}) contains the sub-complex
 \begin{equation}\label{eq:070500}
0\to S^lT^*_aM\ot\tg_a\to S^{l-1}T^*_aM\ot\tg_a\ot\tg_a^*\to
S^{l-2}T^*_aM\ot\tg_a\ot\La^2\tg^*\to\dots
 \end{equation}
and the quotient complex is exact. Indeed, it is the sum by $k$ of
the complexes, each of which is $S^{l-k}(\op{Ann}\tg)$ tensorially
multiplied by the exact complex
 \begin{equation*}
0\to S^k\tg^*_a\ot T_aM/\tg_a\to S^{k-1}\tg^*_a\ot
T_aM/\tg_a\ot\tg_a^*\to S^{k-2}\tg^*_a\ot
T_aM/\tg_a\ot\La^2\tg^*\to\dots
 \end{equation*}
Since complex (\ref{eq:070500}) has nontrivial only
zero-cohomology group, which is isomorphic to
$S^l(\op{Ann}\tg_a)\ot\tg_a$, the same holds for complex
(\ref{hhh2}). Therefore we get ${\frak O}^l=0$ for every $l$,
whence all submanifolds $N$ are transversal.

Similarly if $\Omega$ is a volume form on $M$, the volume
preserving pseudogroup $G=\op{Diff}_\text{loc}(M,\Omega)$ is
involutive: The only non-zero $\d$-cohomology groups are
$H^{0,j}(\g)$. Indeed, $G^1(a_1)=SL(T_aM)$, so
$\g_1=\op{sl}(T_aM)$ and $\g_l=S^lT^*_aM\ot T_aM$ for $l>1$ as for
the general pseudogroup. The homological calculations are similar
and we again get ${\frak O}^l=0$ for every $l$.

Thus for both pseudogroups all submanifolds of codimension $r$ are
locally $G$-equivalent (the pseudogroups act transitively).

\subsection{\hpss Complex pseudogroup}\label{Cpseud}

 \abz
Let $G=\op{Diff}^\CC_{\text{loc}}(M)$ be a pseudogroup of local
holomorphic transformations of a complex manifold $(M,J)$ of
$\dim_\CC M=n$. This group is also involutive, $\g_a^l=S^l_\CC
T_a^*M\ot_\CC T_aM$. Condition (\ref{dim-trans}) reads:
 $$
\tbinom{l+n-1}l\cdot2n\ge\tbinom{l+2n-r-1}l\cdot r.
 $$
This holds true when $r\ge n$, but for $r<n$ it is wrong when
$l>1$. Note though that the for $l=1$ the above inequality holds
for all $0<r<n$. In this case $(\ref{dim-trans})$ is not
sufficient for 1-transversality of all $N$, but it is sufficient
for submanifolds $N$ of general type at the point $a$. For $r\ge
n$ we will show that a generic submanifold $N$ is transversal.

Suppose that $\tg_a$ contains no nontrivial $J$-invariant
subspaces.

Let $T_aM=\tg_a\oplus J\tg_a\oplus\Pi_a$ be a (non-canonical)
decomposition, where $\Pi$ is some $J$-invariant complement to the
first two summands. We get the decomposition ${\frak
h}^l=(J\tg)^*\circ S^{l-1}_\CC T^*M\ot_\CC TM+ \Pi^*\circ_\CC
S^{l-1}_\CC T^*M\ot_\CC TM+S^l_\CC T^*M\ot\tg$. To be more precise
we have two exact three-sequences ($\z^l$ is given by the
diagram):
 $$
  \begin{array}{rcl}
0\to S^l_\CC T^*M\ot\tg\to{\mathfrak h}^l\to\!\!\!&\z^l&\!\!\!\to0\\
&\|&\\
0\to(J\tg)^*\circ S^{l-1}_\CC T^*M\ot
TM/\tg\to\!\!\!&\z^l&\!\!\!\to \Pi^*\circ_\CC S^{l-1}_\CC
T^*M\ot_\CC TM\to0.
  \end{array}
 $$
The Spencer sequences
 $$
\dots\stackrel\d\longrightarrow
g_{i+1}\ot\La^j\tg^*\stackrel\d\longrightarrow
g_i\ot\La^{j+1}\tg^*\stackrel\d\longrightarrow\dots
 $$
with $g_k$ being equal $S^k_\CC T^*M\ot\tg$, $(J\tg)^*\circ
S^{k-1}_\CC T^*M\ot TM/\tg$ or $\Pi^*\circ_\CC S^{k-1}_\CC
T^*M\ot_\CC TM$ has vanishing positive cohomology groups
$H^{i,j+1}=0$, $j\ge 0$. Hence the positive cohomology of
(\ref{hhh2}) vanish too.

Thus Theorem \ref{H-clas} implies the transversality of $N$ with
respect to $G$:

 \begin{prop}\po\label{CC}
The pseudogroup $G$ is $l$-transversal near $a_l=[N]_a^l$ iff\/:
 \begin{enumerate}
\item
 $\tg_a\cap J\tg_a=\{0\}$ or $\tg_a+J\tg_a=T_aM$, when $l=1$,
\item
 $\tg_a\cap J\tg_a=\{0\}$, when $l>1$. \qed
 \end{enumerate}
 \end{prop}

In particular, the complex pseudogroup $G$ acts transitively on
local submanifolds $N\ni a$ of dimension $\op{dim}_\R N\le n$ with
a generic 1-jet at $a$.

On the other hand, a submanifold of dimension $\op{dim}_\R N>n$ is
never transversal. Namely, the intersection $\Pi_a=\tg_a\cap
J\tg_a\ne\{0\}$ and so $N$ possesses an intrinsic geometry.
Investigation of manifolds $N$ equipped with a complex structure
on a distribution $\Pi$ is the subject of Cauchy-Riemann geometry.

The space of covariants ${\frak O}^l_a$ is now non-zero. Thus
there are differential invariants of the action. They are the
curvatures of Cartan-Chern-Moser \cite{CM}. Fixing the curvatures
we get a smaller class $\mathfrak{N}$ of submanifolds, on which
the action will be eventually transitive.

Another important class $\mathfrak{N}$ consists of all complex
submanifolds $N\subset M$ of $\CC$-codimension $r$. This class is
$l$-transversal for every $l$, and so is transitive.
 \vspace{4pt}

\subsection{\hpss Symplectic pseudogroup}\label{S3.4}

 \abz
Consider a symplectic manifold $(M,\oo)$ of dimension $2n$. Let
$G=\op{Symp}(M,\oo)$ be its (pseudo)group of symplectomorphisms.
This pseudogroup is involutive, i.e. $H^{i,j}(\g)=0$ for $i>0$.

Using the identification $TM\stackrel\oo\simeq T^*M$ we write the
symbols $\g^l_a=S^{l+1}T^*_aM\subset S^lT^*_aM\ot T^*_aM$,
understood as homogeneous generating functions (Hamiltonians) of
degree $l+1$.

Condition (\ref{dim-trans}) for $l$-transversality in this case
always holds:
 $$
\tbinom{l+2n}{l+1}\ge\tbinom{l+2n-r-1}{l}\cdot r,
 $$
Thus we get no restrictions on the dimension of a submanifold $N$.

We consider the case of even-dimensional $\tg$, $\dim\tg=2r$ (the
odd-dimensional case is similar). Suppose that $\oo^r$ is
non-degenerate on $\tg$.

Let $TM=\tg\oplus\Pi$ be the direct $\oo$-orthogonal
decomposition, $\Pi=\tg^{\perp_\oo}$. Under the identification
$TM\simeq T^*M$ we have: $\tg\simeq\tg^*$, $\Pi\simeq\Pi^*$. So we
calculate ${\frak h}^1=S^2\tg^*+S^2\Pi^*$. For the prolongations
we have ${\frak h}^l=S^{l+1}\tg^*+S^2\Pi^*\circ_{\text{sym}}
S^{l-1}T^*M$.

Therefore complex (\ref{hhh2}) is the sum of the exact sequence
 $$
0\to S^{l+1}\tg^*\to S^l\tg^*\ot\tg^*\to
S^{l-1}\tg^*\ot\La^2\tg^* \to\dots
 $$
and the tensor product of $S^2\Pi^*$ with the complex
 $$
0\to S^{l-1}T^*M\to S^{l-2}T^*M\ot\tg^*\to
S^{l-3}T^*M\ot\La^2\tg^*\to \dots,
 $$
which (cf.\ (\ref{eq:070500})) has only one nontrivial cohomology
group $H^{l-1,0}\simeq S^{l-1}\Pi^*$.

Thus we again obtain a transitivity result:

 \begin{prop}\po\label{symplpro}
$G$ acts $l$-transversally for all $l\ge0$ (and hence
transitively) near $a_l\in J^l_r(M)$ iff the restriction of $\oo$
to $\tg_a=a_1$ is of maximal rank. \qed
 \end{prop}

The obtained fact is equivalent to a particular case of
Weinstein-Givental theorem (\cite{AG}). To obtain the more general
case we should allow various ranks for the restrictions
$\left.\oo\right|_N$. Then the transversality fails and we get a
1-variant, which is obviously the rank (or dimension of
$\op{Ker}(\left.\oo\right|_N)$). Fixing it we obtain the
transversality for the corresponding equation $\mathfrak{N}$ on
submanifolds.

Finally consider the class $\mathfrak{N}$ of isotropic or
co-isotropic submanifolds. Similar calculations show that $G$ acts
on it $l$-transversally for every $l$.

\subsection{\hpss Contact pseudogroup}\label{cntcp}

 \abz
Consider a contact manifold $(M,\Pi^{2n})$, $\op{dim}M=2n+1$, and
denote by $\nu=TM/\Pi$ the normal. Let $G=\op{Cont}(M,\Pi)$ be the
(pseudo) group of contact transformations. Again the pseudogroup
is involutive. Its Lie algebra consists of contact vector fields
$X_f$, which are determined by generating functions (Hamiltonians)
$f\in C^\infty(M)\ot\nu$.

Necessary condition (\ref{dim-trans}) for $l$-transversality again
holds always:
 $$
\tbinom{l+2n+1}{l+1}\ge\tbinom{l+2n-r}{l}\cdot r,
 $$
Thus no restriction on the dimension of a submanifold $N$ is
imposed.

A choice of a non-zero section of $\nu$ is equivalent to a choice
of a contact form $\a\in C^\infty(\op{Ann}(\Pi)\setminus 0)$,
$\a\we d\a^n\ne0$. Then the Hamiltonian is scalar-valued, $f\in
C^\infty(M)$, and the contact field is uniquely given by
 $$
\a(X_f)=f,\ d\a(\cdot,X_f)=\left.df\right|_{\Pi}.
 $$
In Darboux coordinates $(q,u,p)$, $\a=du-p_idq^i$, we have:
 $$
X_f=\D_{q^i}(f)\p_{p_i}-\p_{p_i}(f)\D_{q^i}+f\p_u, \quad \text{
where }\ \D_{q^i}=\p_{q^i}+p_i\p_u.
 $$

Note that fixing $\a$ is equivalent to the splitting
$T_aM=\Pi_a\oplus\nu_a$, where the first summand is symplectic and
the second is Euclidean 1-dimensional. To describe the symbol
$\g^l_a$, we identify $TM\simeq T^*M$ summand-wise via the
symplectic structure on $\Pi$ and the Euclidean structure on
$\vg$. Then we get:
 $$
\g^l\simeq S^l\nu^*\oplus\sum_{i>0} S^i\Pi^*\ot
S^{l+1-i}\nu^*\simeq S^{l+1}T^*M.
 $$
In fact, order $l$ contact fields $X_f$ are determined by
Hamiltonians $f$ of order $(l+1)$ in all variables except the pure
power of $u$, where the degree is $l$.

Thus the cohomological calculations are quite similar to the
symplectic case and we get:

 \begin{prop}\po
The action of $G$ is $l$-transversal near $a_l\in J^l_r(M)$ for
all $l\ge0$ iff $\tg_a=T_aN$ is transversal to the contact plane
$\Pi_a$ and the induced structure on $\Pi_a^N=\Pi_a\cap T_aN$ from
the canonical conformally-symplectic structure on $\Pi$ is
maximally nondegenerate. This means that
 \begin{itemize}
  \item[-] if $\op{dim}\Pi^N_a=2r$, then
$\left.(d\a)^r\right|_{\Pi^N_a}\ne0$.
  \item[-] if $\op{dim}\Pi^N_a=2r+1$, then
$\op{rank}(\left.d\a\right|_{\Pi^N_a})=2r$. \qed
 \end{itemize}
 \end{prop}
These conditions are equivalent to the claim that through every
point close to $a\in N$ there passes an isotropic submanifold of
dimension no greater than $r$.

As in the symplectic case we note that
$\op{rank}(\left.(d\a)^r\right|_{\Pi^N_a})$ is a 1-variant, fixing
which we get transversality. This is a particular case of the
contact Weinstein-Givental theorem: If two local submanifolds of
codimension $r$ of a contact manifold $(M,\Pi)$ have isomorphic
restrictions of the contact structure $(N,\Pi^N)$, then they have
contactomorphic neighborhoods.

At last, as in \S\,\ref{S3.4}, a particular case says that
restricting to the class $\frak{N}$ of isotropic submanifolds of
fixed dimension, we get transversality of the $G$-action.

\subsection{\hpss Riemannian pseudogroup}

 \abz
Consider at first the isometry pseudogroup of the Euclidean space
$\R^n$. It integrates to the group $G=O(n)\leftthreetimes\R^n$.
The pseudogroup is of finite type and $\g^l_a=0$ for $l\ge2$
(\cite{Ko}). This means that we have plenty of covariants ${\frak
O}^l_a=h^l_a$, $l\ge2$. Thus the transversality is absent for
$l\ge2$, but the action is 1-transversal near each 1-jet $\vp_1$.

Consider now a Riemannian manifold $(M^n,q)$ and let $G$ be the
isometry pseudogroup. If the sectional curvature is constant,
everything is the same is above. But in general case the
pseudogroup $G$ is no longer integrable.

For $l=1,2$ the group $G^l$ is the same as in Euclidean case. But
not every point $\vp_2$ has a prolongation to $G^3$. Indeed,
$G^{1,\text{new}}_a=\rho_{3,1}(G^3_a)\subset G^1_a$ consists of
linear isometries from the orthogonal group $O(T_aM,q)$ preserving
the Riemannian curvature tensor $R_q$.

Moreover, the prolongation-projection method reduces soon the
pseudogroup to the unit element: For a generic Riemannian
structure $q$ the pseudogroups $G^l$ consists of the identity
only, when $l>3$ or $l=3,n>2$.

 \begin{prop}\po
The action of $G$ is not transversal near any $l$-jet for
$l>1$.\qed
 \end{prop}

In fact, various intrinsic and extrinsic curvatures are
$l$-variants. More generally, all differential invariants can be
obtained from the curvatures via Levi-Civita connection operator.
We will obtain transitivity of the action on the equation, which
gives constancy of all these invariants. More generally they give
a solution to the equivalence problem.

\subsection{\hpss Almost complex pseudogroup}\label{S3AC}

 \abz
Consider now an almost complex manifold $(M^{2n},J)$, $J^2=-\1$.
The pseudogroup $G^l$ consists of all $J$-holomorphic $l$-jets:
$J\circ\vp_l=\vp_l\circ J$. It is non-integrable whenever the
almost complex structure $J$ is non-integrable. Let us investigate
this case.

Denote by $N_J\in\op{Hom}_{\bar\CC}(\La^2TM,TM)$ the Nijenhuis
tensor of the structure $J$. This is the obstruction for $J$ to be
integrable (the notation means that it is (2,1)-tensor
$J$-antilinear by each argument).

For $l=1$ we have: $G^1_a=\op{GL}_\CC(T_aM)$, as in the complex
case. The prolongation $G^2$ does not exist over all points of
$G^1$. Using the prolongation-projection method we obtain:
 $$
G^{1,\text{new}}_a=\rho_{2,1}(G^2_a)=\{\Phi\in T^*_aM\ot T_aM
\,|\, J\c\Phi=\Phi\c J,\ N_J\c(\Phi\we\Phi)=\Phi\c N_J \}.
 $$
A symmetric torsion-free connection $\nabla$ on $M$ gives a
decomposition of the 2-jet $\vp_2\in G^2$ into components
$(a,\Phi,\Phi^{(2)})$. The last terms $\Phi^{(2)}\in\Fg(\vp_1)$
for $\vp_1=(a,\Phi)$ are jointly described by the formula:
 $$
\{\Phi^{(2)}\in S^2T^*_aM\ot T_aM,\
J\P^{(2)}(\x,\e)-\P^{(2)}(J\x,\e)=\P\c\nabla_\e
(J)(\x)-\nabla_{\P\e}(J)(\P\x)\},
 $$
Thus $\g_a^2=S^2T^*_aM\ot_\CC T_aM$ as in the complex case, but
for a smaller set of $\vp_1$. The 2-pseudogroup $G^2$ is not
2-integrable in general. Proof of these facts, as well as a
description of the projection $\rho_{l,l-1}:G^l_a\to G^{l-1}_a$
are contained in \cite{Kr1}.

It can be shown, see \cite{Kr2}, that for a generic structure $J$
the set $G^2$ consists of identity for $n>3$, $G^3$ consists of
identity for $n>2$ and $G^4$ is the identity (we ignore the case
$n=1$ corresponding to always integrable $J$). The analysis of
pseudoholomorphic invariants for jets of submanifolds based on the
classification of Nijenhuis tensors (\cite{Kr2}) results in:

 \begin{prop}\po
Let $(M,J)$ be an almost complex manifold with a generic
non-integrable structure $J$ and $n>1$. For $l=1$ the
transversality is described by condition 1 of proposition
\ref{CC}. For $l=2$ no 2-jet is transversal save for the case
$n=2$ and $\dim_\R N=1$. The transversality is absent for $l=3$
and higher. \qed
 \end{prop}

So in the case of generic almost complex pseudogroup we have:
$\g^l=0$ for big $l$, whence plenty of covariants ${\frak O}^l_a$
and no transversality. All the differential invariants here can be
obtained from the Nijenhuis tensor $N_J$ \cite{Kr1}.


\section{\hps Equivalence of differential equations}\label{S4}

 \abz
Let $\pi$ be a vector bundle. A submanifold in $J^k\pi$ can be
identified with a differential equation (actually a system of
equations, but we will just say "equation"). For regularity
purposes we assume it is a subbundle w.r.t.\ all
$\pi_{j,j-1}$-projections, i.e. we have a sequence $\E_j\subset
J^j\pi$ of submanifolds and projections $\pi_{j,j-1}:J^j\pi\to
J^{j-1}\pi$, $j\le k$, forming vector bundles.

Consider two such differential equations $\E\subset J^k\pi$ and
$\E'\subset J^k\pi'$ and two points $x_k$, $x_k'$. There exists a
Lie transformation $\vp:J^\epsilon\pi\to J^\epsilon\pi'$, for
which $\vp^{(k)}(x_k)=x_k'$. So we reduce the problem to the case,
when both equations $\E$ and
$\E'_\vp=\left(\vp^{-1}\right)^{(k)}(\E')$ live in one space
$J^k\pi$. Then we try to identify $\E$ to $\E'_\vp$ by means of a
Lie transformation.

In this section we call Lie pseudogroup $G$ the pseudogroup of Lie
transformations on the jet-spaces $J^k\pi$. Lie-B\"acklund theorem
(\cite{KLV}) states that such a transformation is lifted from a
diffeomorphism of $J^0\pi$ in the case $\op{rank}\pi>1$ (point
transformations) or from a contact diffeomorphism of $J^1\pi$ for
$\op{rank}\pi=1$ (contact transformations). So Lie transformations
are lifted from $J^\epsilon\pi$, where
$\epsilon=\max(0,2-\op{rank}\pi)$.

\subsection{\hpss Formally transitive actions of the Lie
pseudogroup}\label{S4.1.}

 \abz
For the Lie pseudogroup $G$ of point transformations or contact
transformations (depending on $\epsilon$) we provide calculation
of the symbols in Appendix \ref{A2}. This implies (\cite{KL2})
that the dimension of the symbols grow as:
 \begin{equation}
 \label{fo1}
\op{dim}\g^l\sim n_0\cdot\dfrac{l^{n_0-1}}{(n_0-1)!},\qquad
n_0=\dim J^0\pi
 \end{equation}
in the case of point transformations ($\epsilon=0$) and
 \begin{equation}
 \label{fo4}
\op{dim}\g^l\sim \dfrac{l^{n_1-1}}{(n_1-1)!},\qquad n_1=\dim
J^1\pi
 \end{equation}
in the case of contact transformations ($\epsilon=1$). Using the
necessary condition (\ref{dim-trans}) we obtain the following
characterization of equations $\E$ on which the Lie pseudogroup
$G$ acts transitively.

 \begin{theorem}\po
The only transversal equations $\E\subset J^k\pi$ w.r.t.\ the Lie
transformation pseudogroup are the following:
 \begin{enumerate}
  \item
$u'_x=\vp(x,u)$, $x\in\R$, $u\in\R^n$.
  \item
$u'_{x^i}=\vp_i(x,u)$, $i=1,\dots,n$, $x\in\R^n$, $u\in\R$.
  \item
$w'_z=\vp_1(z,\bar z,w,\bar w)$, $w'_{\bar z}=\vp_2(z,\bar
z,w,\bar w)$, $z,w\in\CC$.
  \item
$u'_{x^i}=\vp_i(x,u,u'_{x^{s+1}},\dots,u'_{x^n})$, $1\le i\le
s<n$, $x\in\R^n$, $u\in\R$.
  \item
$u''_{xx}=\vp(x,u,u'_x)$, $x\in\R$, $u\in\R$.  \hfill $\square$
 \end{enumerate}
 \end{theorem}

Let us comment the five cases of the theorem and indicate for
which $\vp$ we actually have the transversality. Let $k$ denote
the order of the equation $\E$, $n$ the dimension of the base of
$\pi$ and $r$ the rank of $\pi$ (dimension of the fiber).

{\bf 1.\/} $n=k=1$. A submanifold $\E\subset J^1(1,r)$ of
codimension $r$ is a determined system of ODEs. Due to our
regularity assumptions they are of main type, so can be written as
in the theorem. By the existence and uniqueness theorem locally
all such systems are equivalent, i.e. the pseudogroup acts
transversally on them. In this case $\E$ is integrable.

{\bf 2.\/} $r=k=1$. Here $\E\subset J^1(n,1)$ of codimension $n$
is diffeomorphically projected by $\pi_0$ to $J^0\pi$. Through
every point $x_0=\pi_0(x_1)$ an $n$-plane $L(x_1)$ passes. Their
collection is the image of the Cartan distribution $\C_\E$ on
$\E$. The obtained rank $n$ distribution on the manifold
$\E^{n+1}$ is generically non-integrable and is either contact or
even-contact. In both cases we get transversality and local
equivalence of all such equations. Thus in this case we take
$\vp_i$ such that $\E$ is maximally non-integrable.

{\bf 3.\/} $k=1$, $n=r=2$. Here the functions should again be
taken generic. Then $\E$ is non-integrable and we obtain
transitivity of the action.

Note that in all the above 3 cases $\pi_{1,0}:\E\to J^0\pi$ is a
diffeomorphism, so that we have the distribution
$d\pi_{1,0}(\C_\E)$ on $J^0\pi$. The described cases correspond to
the known distributions $\Pi$ without moduli: (1) Line field of
$\op{rank}(\Pi)=1$; (2) Contact or even-contact distribution of
$\op{corank}(\Pi)=1$; (3) The Engel distribution of
$\op{rank}(\Pi)=\op{corank}(\Pi)=2$ on a four-dimensional
manifold.

{\bf 4.\/} $k=1$ and $\E\subset J^1(n,1)$ is a submanifold in a
contact manifold. As in \S\ref{cntcp} we see that PDE $\E\subset
J^1\pi$ of $\op{dim}\E=d$ is transversal w.r.t.\ the Lie
pseudogroup at $x_1\in\E$ iff there are no integral manifolds of
the contact structure $\Pi$ of dimension greater than
$\left[\frac{d-1}2\right]$. Note that the induced distribution
$\Pi\cap T\E$ on $\E$ has always integral submanifolds $L$ of
dimension $\left[\frac{d-1}2\right]$. If $\pi_1:L\to L_0$ is a
diffeomorphism, the submanifold has the form $j_1s(L_0)$ for some
section $s$ of the bundle $\pi$. So transversality of $\E$ means
there are no "partial solutions" $s:L_0\to J^0\pi$,
$j_1s(L_0)\subset\E$, of dimension greater than the minimal
possible.

{\bf 5.\/} $k>1$, $n=1$. If $\E_1=\pi_{k,1}(\E)\subset J^1\pi$ is
proper then either there does not exist the prolongation
$\E_1^{(1)}$ or the equation $\E_1$ and hence $\E$ is not
transversal. So we consider $\E_1=J^1\pi$ and then
$\pi_{2,1}:\E_2\to J^1\pi$ is a diffeomorphism. In this case we
have $l$-transversality for every $l$. Note that here $\E$ is
integrable.

The last case corresponds to a known result of S.\ Lie: All the
scalar ordinary differential equations of the second order are
contact equivalent. This case is equivalent to a Legendrian
foliation of the contact 3-manifold $J^1(1,1)$.

In fact, locally all Legendrian foliations of a contact manifold
$J^1(n,1)$ are equivalent, but only for $n=1$ (corresponding to
S.\ Lie's theorem) the corresponding equation $\E$ is generic.
Otherwise, an additional assumption of integrability should be
imposed on $\E$. For such class of equations we get the following
result:

 \begin{prop}\po
The Lie transformations pseudogroup acts transitively on the class
$\mathfrak{N}$ of integrable equations $\E\subset
J^{1+\epsilon}\pi$, such that
$\pi_{1+\epsilon,\epsilon}:\E\stackrel\sim\to J^\epsilon\pi$ is a
diffeomorphism.
 \end{prop}

Indeed, for $\epsilon=1$ we have a Legendrian foliation of
$J^1\pi$ and for $\epsilon=0$ just a foliation of $J^0\pi$.

\subsection{\hpss Formally intransitive actions of the Lie pseudogroup}

 \abz
In all other cases except for the above five, the differential
equations have invariants. Their growth is governed by the
cohomology of covariants $H^{*,*}({\frak O})$.

By Lie-Tresse theorem the invariants of the differential equation
have a finite set of generators. Let us consider some examples.

{\bf 1.\/} Consider a scalar second order ODE $y''=u(x,y,y')$. By
the results of \S\ref{S4.1.} the pseudogroup of contact
transformations acts transitively. For point transformations there
are differential invariants (\cite{Tr2}). The following functions
are the basic differential invariants and the others are obtained
by certain invariant differentiations \cite{C2}.
 $$
I_1=u_{1111},\qquad I_2=u_{xx11}-u_1u_{x11}-4u_{x01}+4u_1u_{01}
-3u_0u_{11}+6u_{00}.
 $$
Here we have denoted the differentiation by $y^{(i)}$ with
subindex $i$, so that we have $u_{x1}=\frac{\p^2u}{\p x\p y'}$
etc.

{\bf 2.\/} Consider the action of $\op{SL}(3)$ on $\R P^2$ by
projective transformations:
 $$
[z_0:z_1:z_2]\mapsto[(Az)_0:(Az)_1:(Az)_2],\quad A\in\op{SL}(3),\
z\in\R^3.
 $$
This action lifts to higher jets $J^k_1(\R P^2)$. It is transitive
near regular orbits for $k\le6$. For $k=6$ it becomes effective.
The first differential invariants $I_7$ occurs for $k=7$. It
equals $(\te_8)^3/(\te_3)^8$, where $\te_3$ and $\te_8$ are some
basic relative differential invariants (\cite{La}; the index $k$
in $\te_k$ refers to the factor under transformations of this
relative invariant). Thus we have an invariant differentiation
$\hat\p=\hat\p_{I_7}$.

The next relative differential invariant is obtained via a bracket
of these $\te_3$ and $\te_8$ (\cite{W}) and using it we can obtain
a new differential invariant $I_8$ of order 8. There will be
exactly one differential invariant of order $k\ge 9$ and each of
them is obtained by the iterated Tresse derivative:
$I_k=\hat\p^{k-8}(I_8)$.

{\bf 3.\/} Monge-Amp\`ere equations with two variables
 $$
\a_0+\a_1u_{xx}+\a_2u_{xy}+\a_3u_{yy}+\a_4\cdot(u_{xx}u_{yy}-u_{xy}^2)=0,
 $$
$\a_i=\a_i(x,y,u,u_x,u_y)$, can be represented geometrically as
effective 2-forms on the contact manifold $J^1(\R^2)$. This gives
a possibility to construct an invariant frame on the equation
($\{e\}$-structure) and so to describe all differential invariants
(\cite{Kr3,KLR}).

{\bf 4.\/} If we consider classification of non-integrable
equations, then their Weyl tensors are differential invariants
w.r.t. Lie pseudogroup of transformations of the jets-space,
\cite{KL1}. Curvatures for geometric structures are particular
cases. Note that prolongation-projection method produces an
integrable equation from a given non-integrable, but the
invariants obtained from the Weyl tensors remain differential
invariants for the new equation.

Let us mention also another related classification problem. Given
a differential equation $\E$ we can consider its pseudogroup of
symmetries, i.e. such transformations from $G$ that map $\E$ to
itself.

These Lie transformations are extrinsic symmetries for the
differential equations. If an equation is not normal (\cite{KLV}),
there may exist also intrinsic symmetries, which cannot be
obtained from the extrinsic ones. Considering a symmetry
pseudogroup $G\subset\op{Sym}(\E)$ of the PDE $\E$, we ask about
equivalence problem for the solutions. The differential invariants
for this problem are important for the integrability of the given
PDE.

For generic $\E$ the symmetry pseudogroup is trivial (this is
possible to show as we did in \S\ref{S3AC} with geometric
structures). But for many important equations $\op{Sym}(\E)$ is
sufficiently big (\cite{KLV}) and the problem is interesting and
non-trivial.

\appendix


\section{\hps Basics from the geometric theory of PDEs}\label{A1}

 \abz
Here we briefly recall some fundamental geometric notions of the
jet-spaces (see \cite{GS1,KLV,Gu,Ly,KL1} for details).

\subsection{\hpss Prolongations and projections}\label{A11}

 \abz
A PDEs system of pure order $k$ is usually represented as a smooth
subbundle $\E\subset J^k\pi$ (\cite{KLV}). This means that
non-regular points are removed and all equations in the system
have pure order $k$. We extend this for different orders.

By a differential equation (system) of maximal order $k$ we mean a
sequence $\E=\{\E_l\}_{-1\le l\le k}$ of submanifolds $\E_l\subset
J^l\pi$ with $\E_{-1}=B$ (base of $\pi$), $\E_0=J^0\pi=E_\pi$ such
that for all $0<l\le k$ the following conditions hold:

 \begin{enumerate}
  \item
$\pi^\E_{l,l-1}:\E_l\to\E_{l-1}$ are smooth fiber bundles.
  \item
The first prolongations $\E_{l-1}^{(1)}$ are smooth subbundles of
$\pi_l$ and $\E_l\subset\E_{l-1}^{(1)}$.
 \end{enumerate}

We remark that in the jets of sections (contrary to the jets of
submanifolds) we have also the projections to the base
$\pi_k:J^k\pi\to B$.

Denote by $\t_x$ the tangent space to the base $B$ of $\pi$ at the
point $x=\pi_k(x_k)$ and by $\nu_{x_0}$ the tangent to the fiber
at $x_0=\pi_{k,0}(x_k)$. Let also $F(x_k)$ be the
$\pi_{k+1,k}$-fiber and $\ups_{x_1}=T_{x_1}\bigl(F(x_0)\bigr)$.

Consider a point $x_k\in\E_k$ with $x_l=\pi_{k,l}(x_k)$ for $l<k$
and $x=x_{-1}$. It determines a symbolic system $g\subset
S\t^*_x\ot\nu_x$ by the formula
 $$
g_l=T_{x_l}\bigl[(\pi^\E_{l,l-1})^{-1}(x_{l-1})\bigr]\subset
S^l\t_x^*\ot\nu_x
 $$
for $l\le k$ and $g_l=g_k^{(l-k)}$ for $l>k$. The conditions above
imply that the symbols $g_l$ form smooth vector bundles over
$\E_l$ and that $g_l\subset g_{l-1}^{(1)}$ for $l\le k$. We call
such collection of subspaces $\{g_k\}$ symbolic systems.

The Spencer $\d$-complex for PDEs system $\E$ at a point
$x_k\in\E_k$ is the Spencer complex for its symbolic system at
this point:
 $$
\cdots\to g_{i+1}\ot\La^{j-1}\t_x^*\stackrel\d\longrightarrow
g_i\ot\La^j\t_x^*\stackrel\d\longrightarrow
g_{i-1}\ot\La^{j+1}\t_x^*\stackrel\d\to\cdots.
 $$
The corresponding $\d$-cohomology is denoted by $H^{i,j}(\E;x_k)$.

We define {\em regular PDEs system\/} of maximal order $k$ as a
submanifold $\E=\E_k\subset J^k\pi$ "cofiltered" by $\E_l$
(property 1 above) and such that the symbolic system and the
Spencer cohomology form graded bundles over it.

Define the {\it Cartan distribution\/} on the space $J^k\pi$ by
the formula: $\mathcal{C}_k(x_k)=(d\pi_{k,k-1})^{-1}L(x_k)$. It
induces the Cartan distribution on $\E_k$:
$\mathcal{C}_{\E_k}=\mathcal{C}_k\cap T\E_k$.

A system of different orders should be investigated for formal
integrability successively by the maximal order $k$. If some
prolongation $\E_k^{(1)}$ is not regular, its projections
$\{\pi_{k+1,l}(\E_{k+1})\}_{l\le k}$, form a new system of maximal
order $k$. Taking the regular part one continues with
prolongations. The process stops in a finite number of steps by
Cartan-Kuranishi theorem on prolongations: There exists a number
$k_0$ such that $\E_k^{(1)}=\E_{k+1}$ for all $k\ge k_0$.

\subsection{\hpss Characteristics}\label{A12}

 \abz
Consider the dual to $g$ system $g^*=\oplus g_k^*$. If $g$ is a
symbolic system, then $g^*$ is an $S\t$-module (as before
$S\t=\oplus S^i\t$ and $\t=\t_x$ with "frozen" $x\in B$) with the
structure given by
 $$
(v\cdot\kappa)p=\kappa(\d_vp),\ v\in S\t,\ \kappa\in g^*,\ p\in g.
 $$
This module, called the {\em symbolic module\/}, is Noetherian and
the Spencer cohomology of $g$ dualizes to the Koszul homology of
$g^*$.

The {\em characteristic ideal\/} is defined by
$I(g)=\op{ann}(g^*)\subset S\t$. The {\em affine characteristic
variety\/} of $g$ (or of $\E$) is the set of $v\in
\t^*\setminus\{0\}$ such that for every $k$ there exists a $w\in
N\setminus\{0\}$ with $v^k\ot w\in g_k$. This is a conical affine
variety. If we consider its complexification and then
projectivization, then we get the {\em characteristic variety\/}
$\op{Char}^\CC(g)\subset P^\CC\t^*$.

Relation of characteristic variety to the characteristic ideal is
given by the formula:
 $$
\op{Char}^\CC(g)=\{p\in P^\CC\t^*\,|\,f(p^k)=0\,\forall f\in I_k,
\forall k\}.
 $$

Note that the dimension of affine characteristic variety equals
the Chevalley dimension of the symbolic module. Recall also that a
sequence of elements $f_1,\dots,f_s\in S\t$ is called regular if
$f_i$ is not a zero divisor in the $S\t$-module
$g^*/(f_1,\dots,f_{i-1})g^*$.

\subsection{\hpss Horizontal differential and generalizations}\label{A13}

 \abz
The horizontal differential $\hat
d:C^\infty(J^k\pi)\to\Omega^1(J^{k+1}\pi)$ is defined by the
properties:
 $$
1.\ \hat df|_{\pi_{k+1}^{-1}(x)}=0,\qquad 2.\ \hat
df|_{j_{k+1}(s)}(x_{k+1})=df|_{j_k(s)}(x_k)
 $$
for any section $s$ of $\pi$ with $j_{k+1}(s)(x)=x_{k+1}$,
$\pi_{k+1,k}(x_{k+1})=x_k$. In local coordinates we can write:
 $$
\hat df=\sum\D_i(f)\,dx^i.
 $$
This can be used as a definition of the total derivative operators
$\D_i$.

Indeed, choosing local coordinates $x^i$ (note the placement of
indices) on the base of $\pi$ and $u^j$ on fibers, we obtain
canonically the coordinates $(x^i,p^j_\z)_{0\le|\z|\le k}$ on
$J^k\pi$, where $p^j_\z([u]_x^k)=\dfrac{\p^{|\z|}u^j}{\p x^\z}$.
Then the operator of total derivative $\D_i:C^\infty(J^k\pi)\to
C^\infty(J^{k+1}\pi)$ has the form:
 $$
\D_i=\p_{x^i}+\sum_{j;\z}p^j_{\z+1_i}\p_{p_\z^j}.
 $$
For a multiindex $\z=(i_1,\dots,i_n)$ we define
$\D_\z=\D_1^{i_1}\cdots\D_n^{i_n}$. If $l=|\z|$ is the length of
the multi-index $\z$, then $\D_\z:C^\infty(J^k\pi)\to
C^\infty(J^{k+l}\pi)$.

If we consider the jet-space $J^k_r({\frak M})$, ${\frak
M}=J^0\pi=E_\pi$, $r=\op{rank}(\pi)$, then $J^k\pi\hookrightarrow
J^k_r({\frak M})$ is an open dense subset. In fact, choosing local
coordinates $(x,u)$ on ${\frak M}$ we identify it locally with
$\pi$ and so the described embedding can be considered as a local
chart. Even though the notions of total derivative and horizontal
differential are not defined on $J^k_r({\frak M})$, we explained
in \S\ref{S13} how to compensate this.


\section{\hps Lie transformations pseudogroup}\label{A2}

 \abz
Consider the pseudogroup of Lie transformations of $M=J^k_r({\frak
M})$. It consists of local diffeomorphisms of the jet-bundle,
preserving the Cartan distribution $\C_k$. This pseudogroup $G$ is
integrable. We will calculate its symbols below.

To simplify we consider the corresponding Lie pseudogroup of
vector fields. A vector field is called an infinitesimal Lie
transformation if its flow is a local Lie transformation. Since we
will work locally, there will be no distinction between
$J^k_r({\frak M})$ and its local chart. Thus we will write
$M=J^k\pi$ for simplicity, where $\pi:E_\pi\to B$ is a vector
bundle.

\subsection{\hpss Lifts of point transformations}\label{A21}

 \abz
In the case $r>1$, denote the projection of the Lie field to
$J^\epsilon\pi=J^0\pi$ by
$X=\sum_ia^i(x,u)\p_{x^i}+\sum_jb^j(x,u)\p_{u^j}$. Then the
prolongation to $J^k\pi$ is
 \begin{equation}\label{X^k}
X^{(k)}=\sum_ia^i(x,u)\D_i^{(k+1)}+\sum_{j;|\z|\le k}
\D_\z(\vp^j)\p_{p^j_\z},
 \end{equation}
where $\vp^j=b^j-\sum^n_{i=1}a^ip^j_i$ are components of the
so-called {\em generating function\/} $\vp=(\vp^1,\dots,\vp^r)$
and $\D_i^{(k+1)}=\p_{x^i}+\sum_{j;|\z|\le
k}p^j_{\z+1_i}\p_{p_\z^j}$ is the operator of total derivative
restricted to $J^k\pi$. Though the coefficients of (\ref{X^k})
depend seemingly on the $(k+1)$-jets, the Lie field is in fact on
$J^k\pi$.

Formula (\ref{X^k}) follows from the claim the Lie field preserves
the co-distribution
 $$
\op{Ann}(\C_k)=\big\langle\oo_\z^j=dp_\z^j-\sum_ip_{\z+1_i}^jdx^i
\,\big|\, 1\le j\le r,\ |\z|<k\big\rangle
 $$
and the formula $d=\sum_i dx^i\ot \D_i^{(k+1)}+ \sum_{j;|\z|\le
k}\oo_\z^j\ot\p_{p_\z^j}$ on $J^k\pi$.

 \begin{prop}\po\label{proppre}
The $l$-symbol $\g^l(x_k)$ of the pseudogroup $G$ at a point
$x_k\in J^k\pi$ admits the splitting $\g^l=\g^l_H\oplus\g^l_V$
depending on a point $x_{k+1}\in F(x_k)$. The horizontal part is
isomorphic to
 $$
\g^l_H(x_k)\simeq
\Bigl[\,S^l\nu_{x_0}^*\oplus\sum_{0<i<k}(S^i\tau_x^*\ot
S^{l-1}\nu_{x_0}^*)\oplus\sum_{i\ge k}(S^i\tau_x^*\ot
S^{k+l-i-1}\nu_{x_0}^*)\Bigr]\ot\tau_x,
 $$
while the vertical (evolutionary) parts is represented as
 $$
\g^l_V(x_k)\simeq \Bigl[\sum_{0\le i<k}(S^i\tau_x^*\ot
S^l\nu_{x_0}^*)\oplus\sum_{i\ge k}(S^i\tau_x^*\ot
S^{k+l-i}\nu_{x_0}^*)\Bigr]\ot\nu_{x_0}.
 $$
 \end{prop}

 \begin{proof}
The space $T_{x_k}J^k\pi$ is decomposed into direct sum of the
horizontal $L(x_{k+1})\subset\C_k(x_k)$ and the vertical
$T_{x_k}^v=\op{Ker}(\pi_k)_*$ components. Thus we have:
 $$
\g^l(x_k)\subset S^lT^*_{x_k}J^k\pi\ot T_{x_k}J^k\pi=
\bigl[S^lT^*_{x_k}J^k\pi\ot L(x_{k+1})\bigr]\oplus
\bigl[S^lT^*_{x_k}J^k\pi\ot T_{x_k}^vJ^k\pi\bigr],
 $$
whence the required splitting. In formula (\ref{X^k}) the
horizontal and vertical components correspond to the first and the
second summands respectively.

Denote by $\mu_a$ the ideal in $C^\infty(B)$ generated by
functions vanishing at $a\in B$, and by $\mu^l_a$ its degree. Let
$\mu_a^l(\mathfrak{Lie})$ be the space of Lie fields vanishing at
$a$ to the order $l$. Then
$\g^l(x_k)=\mu_{x_k}^l(\mathfrak{Lie})/\mu_{x_k}^{l+1}(\mathfrak{Lie})$.

As in the contact and symplectic cases we represent the symbol via
the jets of generating functions. It embeds into the space
$S^lT^*_{x_k}J^k\pi\ot T_{x_k}J^k\pi$ by (\ref{X^k}).

Let us choose a coordinate system such that the point $x_k$
becomes the origin. If $x_k=[s]^k_a$ for some section $s$, this is
achieved by making it the zero-section: $s=\{u^j=0\}$. Then the
condition $X^{(k)}\in\mu^l_{x_k}$ is expressed via the components
of the generating function as follows:
 $$
a^i\in\mu^l_{x_0},\ \p_{x^\z}(a^i)\in\mu^{l-1}_{x_0},\qquad
b^j\in\mu^l_{x_0},\ \p_{x^\z}(b^j)\in\mu^l_{x_0},\qquad
0\le|\z|\le k.
 $$
This yields the claim. Note that the decomposition
$T_{x_0}^*J^0\pi=\tau_x^*\oplus\nu_{x_0}^*$ is induced by the
point $x_1$ and so the representation in the statement is
canonical.
 \end{proof}

\subsection{\hpss Lifts of contact transformations}\label{A22}

 \abz
A Lie transformation for $r=1$ is determined by a contact
transformation $X^{(1)}=X_\vp$ on $J^1\pi$ with a generating
scalar-valued function $\vp=\vp(x^i,u,p_i)$:
 $$
X^{(1)}=\sum_i\Bigl[\D_i^{(1)}(\vp)\p_{p_i}-\p_{p_i}(\vp)\D_i^{(1)}\Bigr]
+\vp\p_u.
 $$
The prolongation of this field to $J^k\pi$ is given by the formula
similar to (\ref{X^k}):
 \begin{equation}\label{X^kk}
X^{(k)}=-\sum_i \p_{p_i}(\vp)\D_i^{(k+1)}+\sum_{|\z|\le
k}\D_\z^{(k)}(\vp)\p_{p_\z}.
 \end{equation}
Again a calculation shows this is a field on $J^k\pi$, coinciding
with $X_\vp$ for $k=1$.

We will need below a decomposition
$T_{x_1}J^1\pi=\tau_x\oplus\nu_{x_0}\oplus\ups_{x_1}$, which is
not canonical. Though the point $x_2$ determines the splitting
$T_{x_1}J^1\pi=L(x_1)\oplus T^v_{x_1}$, the last summand is
further decomposed by a connection in the bundle $\pi_{1,0}$.

 \begin{prop}\po\label{proppro}
The $l$-symbol of the pseudogroup $G$ at a point $x_k\in J^k\pi$
is
 \begin{multline*}
\hspace{-8pt}\g^l(x_k)\simeq \sum_{0\le j\le l}(S^j\nu_{x_0}^*\ot
S^{l+1-j}\ups_{x_1}^*)\oplus S^l\nu_{x_0}^*\oplus\sum_{1\le
i<k;j}(S^i\tau_x^*\ot S^j\nu_{x_0}^*\ot S^{l-j}\ups_{x_1}^*)\\
\qquad\oplus \sum_{k\le i;j}(S^i\tau_x^*\ot S^j\nu_{x_0}^*\ot
S^{k+l-i-j}\ups_{x_1}^*).
 \end{multline*}
 \end{prop}

 \begin{proof}
As in proposition \ref{proppre}, due to (\ref{X^kk}), in the
coordinate system $(x^i,u)$ such that $p_\z(x_k)=0$ for $|\z|\le
k$, the condition $X^{(k)}\in\mu^l_{x_k}(\mathfrak{Lie})$ is
equivalent to:
 $$
\vp\in\mu^l_{x_1},\ \p_{p_i}(\vp)\in\mu^l_{x_1},\
\p_{x^\z}(\vp)\in\mu^l_{x_1},\qquad 0\le|\z|\le k.
 $$
The claim follows.
 \end{proof}


\end{document}